\documentclass[letterpaper, 10 pt, conference]{ieeeconf}
\IEEEoverridecommandlockouts
\usepackage{amsmath,amssymb,amsthm}
\usepackage{graphicx}
\usepackage[protrusion=false]{microtype}
\usepackage{booktabs}
\usepackage{array}
\usepackage{placeins}
\usepackage{cite}
\usepackage{algorithm}
\usepackage{algpseudocode}
\algrenewcommand\algorithmicrequire{\textbf{Input:}}
\algrenewcommand\algorithmiccomment[1]{\hfill{\itshape #1}}
\algrenewcommand\alglinenumber[1]{\footnotesize #1}
\algnewcommand\Parameters{\item[\textbf{Parameters:}]}
\algtext*{EndFor}
\makeatletter
\renewcommand\floatc@ruled[2]{{\small\bfseries #1}\ \small #2\par}
\let\@origmakecaption\@makecaption
\long\def\@makecaption#1#2{\ifx\@captype\@IEEEtablestring\begin{center}{\footnotesize #1:~\scshape #2}\end{center}\@IEEEtablecaptionsepspace\else\@origmakecaption{#1}{#2}\fi}
\makeatother
\usepackage{xcolor}
\usepackage{changepage}
\newcommand{\repolink}{\mbox{\footnotesize\url{https://github.com/marginallyunstable/klpi4alqr}}}
\usepackage[hidelinks]{hyperref}

\theoremstyle{plain}
\newtheorem{theorem}{Theorem}
\newtheorem{proposition}{Proposition}
\newtheorem{corollary}{Corollary}

\newtheorem{assumption}{Assumption}
\newtheorem{applemma}{Lemma}
\counterwithin*{applemma}{section}

\newtheorem{fact}{Fact}

\theoremstyle{definition}
\newtheorem{remark}{Remark}

\newcommand{\Tr}{\operatorname{Tr}}
\newcommand{\KL}{D_{\mathrm{KL}}}
\newcommand{\Nor}{\mathcal N}
\newcommand{\EE}{\mathbb E}
\newcommand{\PP}{\mathbb P}
\newcommand{\cS}{\mathcal S}
\newcommand{\Qf}{\mathcal Q}
\newcommand{\usig}{\underline\sigma}
\newcommand{\SNR}{\mathrm{SNR}}
\newcommand{\sa}{\bar s_{\mathrm a}}
\newcommand{\Wb}{\overline W}
\newcommand{\pitem}[1]{\par\noindent#1\ }
\newcommand{\runin}[1]{\par\smallskip\noindent\textbf{#1:}\ }

\renewcommand{\baselinestretch}{0.978}
\newlength{\EqSkipAbove}
\newlength{\EqSkipBelow}
\newlength{\EqSkipShortAbove}
\newlength{\EqSkipShortBelow}
\makeatletter
\newcommand{\ApplyEqSkips}{%
  \setlength{\abovedisplayskip}{\EqSkipAbove}%
  \setlength{\belowdisplayskip}{\EqSkipBelow}%
  \setlength{\abovedisplayshortskip}{\EqSkipShortAbove}%
  \setlength{\belowdisplayshortskip}{\EqSkipShortBelow}}
\AtBeginDocument{\ApplyEqSkips\g@addto@macro\normalsize{\ApplyEqSkips}}
\makeatother
\title{\LARGE \bf KL-Regularized Policy Iteration for Adaptive LQR:\\ Gain Updates and Probing from One Policy Improvement}

\author{Ajinkya Bhole$^{1,2}$, Mohammad Mahmoudi Filabadi$^{1,2}$, Guillaume Crevecoeur$^{1,2}$, and Tom Lefebvre$^{1,2}$%
\thanks{$^{1}$Department of Electromechanical, Systems and Metal Engineering, Ghent University, Ghent, Belgium.}%
\thanks{$^{2}$Core lab MIRO, Flanders Make, Belgium.}}

\begin{document}
\maketitle
\thispagestyle{empty}
\pagestyle{empty}

\begin{abstract}
Policy-gradient methods for adaptive linear quadratic control improve the feedback gain of an unknown plant online by gradient-type steps computed from closed-loop data.
A recent framework of such methods, policy-gradient adaptive control (PGAC)~\cite{zhao2025pgac}, takes one gradient-type step per sample, which its indirect form computes from an estimated model and its direct form directly from the data.
Because data generated by a fixed linear feedback cannot identify the plant, PGAC adds a probing perturbation to the input and leaves its covariance as a design parameter.
We study the adaptive LQR problem through a KL-regularized formulation, in which the probing is the randomness of a Gaussian policy kept close to a Gaussian reference policy, so that the cost of the probing is part of the objective.
One policy improvement of this formulation gives a gain update together with a probing covariance.
We show that, for suitable references centred at the current gain, this gain update is exactly PGAC's Gauss--Newton or natural-gradient update with any step size below a bound.
We repeat this policy improvement, each time against a tempered copy of the current policy, and refer to this scheme as KL-regularized policy iteration (KLPI).
On a known plant and with tempering, the cost of KLPI's gains converges linearly to the optimal cost, and the probing covariance keeps a floor proportional to the regularization weight.
On an estimated plant, under conditions like PGAC's (small updates and sufficiently informative data), the closed loop stays stable, the cost of the gains converges linearly up to a bias, and the probing covariance again keeps this floor.
In simulations at equal probing cost, the shape of KLPI's probing covariance, which PGAC can adopt, lowers the extra cost caused by suboptimal gains on nearly all tested plants, whereas KLPI's use of past iterations in its gain update changes this cost little.
\end{abstract}

\section{Introduction}\label{sec:intro}

The linear quadratic regulator (LQR) asks for the state-feedback gain $K$ that minimizes the infinite-horizon average quadratic cost $C(K)$ of a linear plant driven by Gaussian process noise.
When the plant is known, the optimal gain $K^\star$ follows from a Riccati equation, or from Hewer's policy iteration~\cite{hewer1971iterative}, which alternates evaluating the current gain (policy evaluation) with choosing the input that minimizes its state-action value (policy improvement).
Policy-gradient methods search over gains directly.
The cost $C(K)$ is not convex in $K$, but Fazel et al.~\cite{fazel2018global} showed that it is gradient dominated, that is, the gap $C(K)-C(K^\star)$ is at most a multiple of the squared norm of the gradient.
With suitable step sizes, the gradient, natural-gradient and Gauss--Newton iterations therefore converge linearly from any stabilizing gain.
The Gauss--Newton step rescales the gradient on the left by the inverse Hessian of the state-action value with respect to the input and on the right by the inverse stationary state covariance, and with step size one half it is Hewer's step~\cite{bu2019lqr}.

In adaptive LQR the plant is unknown and the controller improves its gain from closed-loop data, which raises two difficulties.
Firstly, the switching between gains that each stabilize the plant can still destabilize the closed loop.
The state stays bounded when the gains are uniformly strongly stable and consecutive gains stay close, a property used in~\cite{cohen2019learning} that we call sequential stability in this paper.
A gradient method secures it by taking small steps.
Secondly, the data generated by a fixed linear feedback $u=Kx$ cannot identify the plant, because several plants produce the same closed loop on the subspace where these data lie~\cite{simchowitz2020naive}.
The controller must therefore add to its input a probing perturbation, that is, a component that is not a linear function of the state.

Policy-gradient adaptive control (PGAC)~\cite{zhao2025pgac} is a framework of policy-gradient methods for adaptive LQR that addresses both difficulties.
At every sample, a PGAC method applies the current gain plus a probing perturbation and then takes one gradient, natural-gradient or Gauss--Newton step on the gain.
The indirect form computes this step from a recursive least-squares estimate of the plant, and the direct form computes it from the sample covariance of the data.
Computing a step requires solving only Lyapunov equations, whereas recomputing the optimal gain of the changing plant estimate at every sample would require solving a Riccati equation.
Zhao et al.~\cite{zhao2025pgac} also cite the robustness of small updates.
Zhao et al.~\cite{zhao2025pgac} prove that, if the step size is small and the data have a large enough signal-to-noise ratio, the switching gains keep the closed loop stable and the cost of the gain converges linearly up to a bias that shrinks as this ratio grows.
In their analysis the probing matters only through the data it produces, and its cost is not part of the cost that the gain update decreases.
The covariance of the probing perturbation is a design parameter, constant and isotropic in the experiments of~\cite{zhao2025pgac}, and both PGAC and the earlier data-enabled policy optimization (DeePO)~\cite{zhao2025deepo} leave its choice to the designer.

This paper studies adaptive LQR through the lens of KL regularization, which is common in reinforcement learning~\cite{geist2019regularized} and in probabilistic approaches to control~\cite{lefebvre2024probabilistic}.
A KL-regularized problem minimizes the expected cost plus a penalty equal to a positive regularization weight times the KL divergence of the policy from a reference policy.
Two of its properties matter for adaptive LQR.
Firstly, when the reference is a Gaussian policy, the optimal policy is also randomized, as for the LQR with an entropy bonus~\cite{guo2023fast}, so the problem \emph{produces the probing} and \emph{counts its cost in the objective}.
Secondly, when the reference is the current policy, the KL penalty keeps the new policy close to it, which gives the small steps that sequential stability needs.
Just as PGAC takes one gradient-type step per sample, we take one policy improvement per sample instead of solving the regularized problem.
Repeating this policy improvement with the current policy as reference and a fixed regularization weight makes the probing die out, because each repetition adds to the precision (inverse covariance) of the policy.
Tempering the reference, that is, raising the density of the current policy to a power $\alpha<1$ and renormalizing, gives a policy improvement that penalizes the KL divergence to the current policy and rewards the entropy~\cite{kozuno2019theoretical,vieillard2020leverage,mahmoudi2026kl}.
For the LQR, tempering keeps the policy covariance bounded away from zero when the regularization weight is bounded away from zero.
For brevity, we refer to this scheme, in which the policy improvement is repeated against the tempered current policy, as KL-regularized policy iteration (KLPI).

\runin{Related work}
KLPI applies to adaptive LQR the KL- and entropy-regularized policy iteration of reinforcement learning~\cite{kozuno2019theoretical,vieillard2020leverage}.
This policy iteration is a policy mirror-descent method~\cite{vieillard2020leverage}.
Without tempering and with a fixed regularization weight, the repeated policy improvement reweights the initial policy exponentially in the summed state-action values of earlier policies, as in POLITEX~\cite{abbasi2019politex} for finite Markov decision processes.
The model-free LQR method MFLQ~\cite{abbasi2019model} instead plays the gain that minimizes its sum of past estimated state-action values and probes with scheduled random inputs.
In trajectory optimization, Mahmoudi Filabadi et al.~\cite{mahmoudi2026kl} regularize the backward pass of differential dynamic programming with an entropy bonus and a KL penalty to the policy of the previous iteration.
Over the iterations, their Gaussian policies accumulate the input gradients and Hessians of the local state-action values with geometrically decaying weights, and the authors interpret this accumulation as momentum.
Closest to this paper, Szpruch et al.~\cite{szpruch2024optimal} probe only through entropy-regularized Gaussian policies when controlling an unknown continuous-time linear-quadratic system over finite-horizon episodes, and one of their algorithms penalizes the KL divergence to the previous episode's policy.
In this paper, every policy improvement evaluates the current gain on a least-squares plant estimate updated at every sample of a single trajectory, and our stability analysis follows that of~\cite{zhao2025pgac}.

\runin{Contributions}
(i) We show that one policy improvement of the KL-regularized LQR returns a gain update and a probing covariance, and that for suitable references the gain update is PGAC's Gauss--Newton or natural-gradient step with any step size below a bound (Propositions~\ref{prop:soft} and~\ref{prop:corr}).
(ii) With a known plant and tempering, we prove that the cost of KLPI's gains converges linearly and that its probing covariance keeps a floor proportional to the regularization weight (Theorem~\ref{thm:mb}).
(iii) On an estimated plant, under conditions like those in PGAC's analysis, we prove that KLPI keeps the closed loop stable and that its cost converges linearly up to a bias (Theorem~\ref{thm:adaptive}), and that its probing covariance keeps a floor proportional to the regularization weight (Corollary~\ref{cor:floor}).
(iv) In simulations at equal probing cost, the shape of KLPI's probing covariance, which PGAC can adopt, lowers the extra cost caused by suboptimal gains on nearly all tested plants, while KLPI's use of past iterations in its gain update changes this cost little (Section~\ref{sec:cmp}).

\runin{Organization}
Section~\ref{sec:setting} recalls the LQR, its policy-gradient updates and PGAC.
Section~\ref{sec:kl} introduces the KL-regularized LQR problem and its policy improvement.
Section~\ref{sec:alg} introduces KLPI, analyses it when the plant is known and compares it with PGAC's Gauss--Newton iterations.
Section~\ref{sec:adaptive} analyses KLPI on an estimated plant.
Section~\ref{sec:cmp} compares KLPI and PGAC in simulations, and Section~\ref{sec:concl} concludes.

\runin{Notation}
$\|\cdot\|$ is the spectral norm, $\usig(\cdot)$ the smallest singular value, $\lambda_{\min}$, $\lambda_{\max}$ the extreme eigenvalues of a symmetric matrix, and $M\succeq N$ means that $M-N$ is positive semidefinite.
$\Nor(m,\Psi)$ is the Gaussian distribution with mean $m$ and covariance $\Psi$.
Capitals $X_t$, $U_t$, $W_t$, $\Omega_t$ denote the closed-loop state, input, process noise and probing perturbation, lowercase $x$, $u$ denote values, and $Z_t:=[U_t;X_t]$.
$\KL(\mu'\|\mu)$ is the KL divergence of a distribution $\mu'$ from a distribution $\mu$, and $h(\mu)$ is the differential entropy of $\mu$.
For stationary policies $\pi$ and $\rho$, we write $\kappa^\pi_\rho(x):=\KL(\pi(\cdot\mid x)\|\rho(\cdot\mid x))$ for the KL divergence of $\pi$ from $\rho$ in the state $x$.

\section{Policy-gradient methods for LQR \& PGAC}\label{sec:setting}

This section sets up the plant and its cost, then states the policy-gradient updates that PGAC takes and its condition on the data, which Section~\ref{sec:adaptive} also adopts.

\runin{Plant and cost}
We consider the plant
\begin{equation}\label{eq:plant}
X_{t+1}=AX_t+BU_t+W_t,\qquad W_t\sim\Nor(0,\Sigma_W),
\end{equation}
with $X_t\in\mathbb R^n$, $U_t\in\mathbb R^m$, a process-noise covariance $\Sigma_W\succ0$, $(W_t)$ independent over time and of the initial state $X_0$, which has a finite second moment, $(A,B)$ controllable, and the stage cost $c(x,u):=x^\top Qx+u^\top Ru$ with $Q,R\succ0$.
A stationary policy $\pi$ assigns to each state $x$ a distribution $\pi(\cdot\mid x)$ of the input, and $\PP^\pi$ denotes the law of the closed loop~\eqref{eq:plant} with $U_t\sim\pi(\cdot\mid X_t)$.
Its average cost is
\begin{equation}\label{eq:avgcost}
C(\pi):=\limsup_{T\to\infty}\tfrac1T\,\EE_{\PP^\pi}\Big[\textstyle\sum_{t<T}c(X_t,U_t)\Big].
\end{equation}
A gain $K$ defines the deterministic policy $u=Kx$ and the closed-loop matrix $A_K:=A+BK$, and $\cS$ is the set of gains for which $A_K$ is Schur.
Adding a Gaussian probing perturbation of covariance $\Psi$, drawn independently of the past, gives the linear-Gaussian policy $\Nor(Kx,\Psi)$.

\runin{Matrices of a gain}
For $K\in\cS$, the cost matrix $P_K$ solves $P_K=Q+K^\top RK+A_K^\top P_KA_K$, the stationary state covariance solves $\Sigma_K=\Sigma_W+A_K\Sigma_KA_K^\top$, so that $\Sigma_K\succeq\Sigma_W$, and $C(K)=\Tr(P_K\Sigma_W)=\Tr((Q+K^\top RK)\Sigma_K)$.
For $K\notin\cS$, $C(K)=\infty$.
The state-action value of $K$ is the cost of applying $u$ once in state $x$ and $Kx$ afterwards, relative to the average cost,
\begin{equation}\label{eq:Q}
\resizebox{0.91\columnwidth}{!}{$\begin{aligned}
\Qf_K(x,u)&:=c(x,u)-C(K)+\EE_{X'\sim\Nor(Ax+Bu,\Sigma_W)}\big[X'^\top P_KX'\big]\\
&\phantom{:}=u^\top S_Ku+2u^\top N_Kx+x^\top(Q+A^\top P_KA)x,
\end{aligned}$}
\end{equation}
with the input Hessian $S_K:=R+B^\top P_KB$ (half the Hessian of $\Qf_K$ in $u$) and $N_K:=B^\top P_KA$.
Its gradient in $u$ at $u=Kx$ is $2E_Kx$, with the slope $E_K:=S_KK+N_K$, which vanishes exactly at $K^\star$.
We write $P^\star$, $S^\star$, $\Sigma^\star$, $C^\star$ for the objects of $K^\star$, and $\gamma^\star:=\|\Sigma^\star\|/\usig(\Sigma_W)\ge1$.
Since $P_K\succeq Q$, every input Hessian satisfies $S_K\succeq \underline{s}I$ with $\underline{s}:=\lambda_{\min}(R+B^\top QB)$.

\runin{Policy-gradient updates of the gain}
The gradient of the cost is $\nabla C(K)=2E_K\Sigma_K$~\cite{fazel2018global}.
With a step size $\eta>0$, the natural-gradient and Gauss--Newton steps of~\cite{fazel2018global} are
\begin{equation}\label{eq:steps}
K-2\eta E_K,\qquad K-2\eta S_K^{-1}E_K.
\end{equation}
Both steps of~\eqref{eq:steps} are preconditioned steps $K-H^{-1}E_K$, in which a symmetric positive definite matrix $H$ (the preconditioner) sets the size and direction of the step, with $H=\tfrac{1}{2\eta}I$ and $H=\tfrac{S_K}{2\eta}$ respectively.
For $\eta=\frac12$ the Gauss--Newton step minimizes~\eqref{eq:Q} over $u$, which is Hewer's step~\cite{hewer1971iterative,bu2019lqr}.

\runin{Adaptive loop and PGAC}
In the adaptive problem the controller knows $Q$, $R$ and the state but not $(A,B)$, and it starts from $t_0$ offline samples and a gain $K_{t_0}\in\cS$.
In what follows, PGAC refers to its indirect form.
At time $t$ PGAC applies $U_t=K_tX_t+\Omega_t$ with probing perturbation $\Omega_t$, observes $X_{t+1}$, updates the least-squares estimate $\hat\Theta_{t+1}=[\hat B_{t+1}\ \hat A_{t+1}]$ of $\Theta:=[B\ A]$, and takes the step~\eqref{eq:steps} with $S_{K_t}$, $E_{K_t}$ computed on this plant estimate~\cite{zhao2025pgac}.
The data enter through $\Phi_t:=\frac1t\sum_{i<t}Z_iZ_i^\top$ and $\Wb_t:=\frac1t\sum_{i<t}W_iZ_i^\top$, so that $\hat\Theta_t-\Theta=\Wb_t\Phi_t^{-1}$, and PGAC measures how informative they are by the signal-to-noise ratio $\SNR_t:=\usig(\Phi_t)/\|\Wb_t\|$.

\begin{assumption}[Invertible data~\cite{zhao2025pgac}]\label{ass:pe}
For all $t\ge t_0$, $\Phi_t$ is invertible.
\end{assumption}

Then $\|\hat\Theta_t-\Theta\|\le1/\SNR_t$~\cite[Lemma~1]{zhao2025pgac}.
Suppose that $\SNR_t$ exceeds a constant for all $t$, which we call the data condition, and that the step size is small.
For the gradient, natural-gradient and Gauss--Newton steps,~\cite[Theorems~1, 3 and~4]{zhao2025pgac} then state that the gains are sequentially stable, that the state is bounded by a decaying transient plus a multiple of the largest disturbance, and that $C(K_t)-C^\star$ contracts linearly up to a bias equal to a discounted sum of $1/\SNR_i$ over $t_0\le i\le t$.

\section{One KL-regularized policy improvement}\label{sec:kl}

This section introduces the KL-regularized LQR (P), then computes one policy improvement of (P) and shows that it returns a gain update together with a probing covariance (Proposition~\ref{prop:soft}).
It then computes what the probing costs (Proposition~\ref{prop:price}) and identifies the references for which the gain update is one of PGAC's steps (Proposition~\ref{prop:corr}).

\runin{The regularized problem}
Let $\rho$ be a reference policy and $\tau>0$ a regularization weight.
The KL-regularized LQR is
\begin{equation}\label{eq:P}
\begin{aligned}
&\text{(P)}\quad \min_\pi\ J_\rho(\pi),\\
&J_\rho(\pi):=\limsup_{T\to\infty}\tfrac1T\,\EE_{\PP^\pi}\Big[\textstyle\sum_{t<T}\big(c(X_t,U_t)+\tau\,\kappa^\pi_\rho(X_t)\big)\Big],
\end{aligned}
\end{equation}
over stationary policies $\pi$.
When the reference has a positive definite covariance, every deterministic policy has an infinite penalty, so (P) selects a randomized policy.
We take Gaussian references of the form of PGAC's policy, $\rho(\cdot\mid x)=\Nor(K_\rho x,\Psi_\rho)$ with $\Psi_\rho\succ0$, and write $D:=\frac\tau2\Psi_\rho^{-1}$ for the scaled precision (inverse covariance) of the reference.

\runin{Policy improvement}
Policy iteration for (P) first evaluates the current policy $\pi$, which yields a relative state value $V_\pi$ (defined up to an additive constant) and the state-action value $\Qf_\pi(x,u):=c(x,u)+\EE_{X'\sim\Nor(Ax+Bu,\Sigma_W)}[V_\pi(X')]$.
It then improves $\pi$, that is, it replaces $\pi$ by
\begin{equation}\label{eq:softpi}
\pi'(\cdot\mid x)\in\arg\min_{\mu}\EE_{U\sim\mu}\big[\Qf_\pi(x,U)\big]+\tau\KL\big(\mu\,\|\,\rho(\cdot\mid x)\big).
\end{equation}
Its minimizer is the reference reweighted by $e^{-\Qf_\pi(x,\cdot)/\tau}$ and normalized (Lemma~\ref{lem:gibbs}).

\begin{proposition}[Policy improvement]\label{prop:soft}
Let $K\in\cS$, $\Psi\succ0$, $\pi=\Nor(Kx,\Psi)$ and $\rho=\Nor(K_\rho x,\Psi_\rho)$ with $\Psi_\rho\succ0$.
Let $P^\rho_K$ solve $P=Q+K^\top RK+(K-K_\rho)^\top D(K-K_\rho)+A_K^\top PA_K$, and let $S^\rho_K:=R+B^\top P^\rho_KB$ and $E^\rho_K:=S^\rho_KK+B^\top P^\rho_KA+D(K-K_\rho)$.
\pitem{(a)} $V_\pi(x)=x^\top P^\rho_Kx$ and $J_\rho(\pi)=\Tr(P^\rho_K\Sigma_W)+\Tr(S^\rho_K\Psi)+\tau\KL(\Nor(0,\Psi)\|\Nor(0,\Psi_\rho))$.
\pitem{(b)} The policy improvement~\eqref{eq:softpi} returns $\pi'=\Nor(K'x,\Psi')$ with
\begin{equation}\label{eq:softgen}
K'=K-H^{-1}E^\rho_K,\quad \Psi'=\tfrac\tau2H^{-1},\quad H:=S^\rho_K+D.
\end{equation}
\pitem{(c)} $K'\in\cS$ and $J_\rho(\pi')\le J_\rho(\pi)$.
\pitem{(d)} If $K_\rho=K$, then $P^\rho_K=P_K$, $\Qf_\pi=\Qf_K+C(K)$, and~\eqref{eq:softgen} reads
\begin{equation}\label{eq:softclosed}
K'=K-H^{-1}E_K,\quad \Psi'=\tfrac\tau2H^{-1},\quad H:=S_K+D,
\end{equation}
with $P_{K'}\preceq P_K-E_K^\top H^{-1}E_K$.
\end{proposition}
\begin{proof}See Appendix~\ref{app:kl}.\end{proof}

The gain update~\eqref{eq:softgen} is Hewer's step for the LQR whose stage cost $c$ is augmented by $(u-K_\rho x)^\top D(u-K_\rho x)$, a penalty on the distance of the input from the mean of the reference.
We call a reference with $K_\rho=K$, as in Proposition~\ref{prop:soft}(d), centred at $K$.
With a centred reference and a fixed regularization weight, the KL penalty pulls the new gain towards the current one, and the pull grows as the reference covariance shrinks.

\begin{proposition}[Cost of the probing]\label{prop:price}
For $K\in\cS$ and $\Psi\succeq0$, the average cost~\eqref{eq:avgcost} of $\Nor(Kx,\Psi)$ is $C(K,\Psi)=C(K)+\Tr(S_K\Psi)$.
For $\Psi\succ0$, $\Tr(S_K\Psi)\ge m(\det S_K\det\Psi)^{1/m}$, with equality if and only if $\Psi$ is proportional to $S_K^{-1}$.
\end{proposition}
\begin{proof}See Appendix~\ref{app:kl}.\end{proof}

Proposition~\ref{prop:price} separates two properties of a probing covariance $\Psi$,
\begin{equation}\label{eq:amountshape}
\text{amount }\ \Tr(S_K\Psi),\qquad\text{shape }\ \{\lambda\Psi:\lambda>0\}.
\end{equation}
The amount of $\Psi$ is the cost that the probing adds to the average cost at every sample, and the shape of $\Psi$ is $\Psi$ up to a positive factor.
Among covariances with a given determinant, and hence a given entropy of $\Nor(0,\Psi)$, the amount is smallest for the shape $S_K^{-1}$.
Proposition~\ref{prop:corr} identifies the centred references for which the gain update of the policy improvement is exactly PGAC's natural-gradient or Gauss--Newton step~\eqref{eq:steps}.

\begin{proposition}[Correspondence]\label{prop:corr}
Let $K\in\cS$ and $\tau>0$.
\pitem{(i)} Every symmetric $H\succ S_K$ is the preconditioner~\eqref{eq:softclosed} of exactly one centred reference, the one with $\Psi_\rho=\frac\tau2(H-S_K)^{-1}$.
\pitem{(ii)} The Gauss--Newton step with $0<\eta<\frac12$ is the gain update of the policy improvement against $\Nor\big(Kx,\frac{\eta\tau}{1-2\eta}S_K^{-1}\big)$, whose probing covariance is $\Psi'=\eta\tau S_K^{-1}$.
\pitem{(iii)} The natural-gradient step with $0<\eta<\tfrac{1}{2\lambda_{\max}(S_K)}$ is the gain update of the policy improvement against the reference $\Nor\big(Kx,\tau(\eta^{-1}I-2S_K)^{-1}\big)$, whose probing covariance is $\Psi'=\eta\tau I$.
\pitem{(iv)} Hewer's step is the limit of~\eqref{eq:softclosed} as $D\to0$, that is, as the reference covariance grows without bound.
If $E_K\neq0$, no Gaussian reference centred at $K$ gives the Gauss--Newton step with $\eta\ge\frac12$.
\end{proposition}
\begin{proof}See Appendix~\ref{app:kl}.\end{proof}

Thus PGAC's natural-gradient step with $\eta<\tfrac{1}{2\lambda_{\max}(S_K)}$, followed by an isotropic probing perturbation of variance $\sigma^2$, is the policy improvement against the centred reference of Proposition~\ref{prop:corr}(iii) with $\tau=\sigma^2/\eta$.
When $E_K$ has full row rank, PGAC's Gauss--Newton step with its probing perturbation is a policy improvement only if the probing covariance is $\eta\tau S_K^{-1}$, and constant isotropic probing has this form only when $S_K$ is a multiple of the identity.

\section{KL-regularized policy iteration with a tempered reference}\label{sec:alg}

This section introduces KLPI.
We first show that, when the plant is known and the regularization weight is constant, no iteration of KLPI increases a fixed objective.
We then prove that the cost of its gains converges, and finally compare KLPI with PGAC.

\runin{The tempered reference}
Taking the current policy itself as the next reference makes the probing die out (Theorem~\ref{thm:mb}(b)), and tempering the policy first prevents this.
To temper a Gaussian policy, we raise its density to a power $\alpha_t>0$ and normalize it, which keeps the mean and divides the covariance by $\alpha_t$.
For $\pi_t=\Nor(K_tx,\Psi_t)$, the tempered reference $\rho_t\propto\pi_t^{\alpha_t}$ is therefore $\Nor(K_tx,\Psi_t/\alpha_t)$, which is centred at $K_t$, so Proposition~\ref{prop:soft}(d) applies.

\begin{proposition}[One iteration of KLPI]\label{prop:alg}
Let $\alpha\in(0,1]$, $\tau_t>0$, $\pi_t=\Nor(K_tx,\Psi_t)$ with $\Psi_t=\frac{\tau_t}2H_t^{-1}$, $H_t\succ0$, and let $S$, $E$ be the input Hessian and slope of $K_t$ on a model that $K_t$ stabilizes.
The policy improvement with regularization weight $\tau_{t+1}$ against $\rho_t\propto\pi_t^{\alpha_t}$, $\alpha_t:=\alpha\tau_t/\tau_{t+1}$, is $\pi_{t+1}=\Nor(K_{t+1}x,\Psi_{t+1})$ with
\begin{equation}\label{eq:alg}
\begin{gathered}
H_{t+1}=\alpha H_t+S,\qquad K_{t+1}=K_t-H_{t+1}^{-1}E,\\
\Psi_{t+1}=\tfrac{\tau_{t+1}}2H_{t+1}^{-1}.
\end{gathered}
\end{equation}
\end{proposition}
\begin{proof}The reference $\rho_t=\Nor(K_tx,\Psi_t/\alpha_t)$ has the scaled precision $D=\frac{\tau_{t+1}}2\alpha_t\Psi_t^{-1}=\alpha H_t$, and Proposition~\ref{prop:soft}(d), applied to the model on which $S$ and $E$ are computed, gives~\eqref{eq:alg}.\end{proof}

With $\alpha=1$ and a constant regularization weight, $\alpha_t=1$, so the reference is the current policy itself, and we call this case untempered.
The regularization weight does not enter the recursion of $(K_t,H_t)$. It only converts $H_t$ into the probing covariance.
Unrolling the recursion gives $(1-\alpha)H_{t+1}=\sum_{k\le t}(1-\alpha)\alpha^{t-k}S_{K_k}+\alpha^{t+1}(1-\alpha)H_0=:\bar S_t$, which for $\alpha<1$ is a weighted average of the past input Hessians and of $(1-\alpha)H_0$.
The gain update $K_{t+1}=K_t-(1-\alpha)\bar S_t^{-1}E$ is thus a Gauss--Newton step with step size $\eta=\frac{1-\alpha}2$ and an averaged input Hessian.

\runin{A fixed objective}
Because the reference changes at every iteration, each policy improvement decreases a different regularized objective.
With a known plant and a constant regularization weight $\tau$, KLPI nevertheless does not increase one fixed objective, defined for $\pi=\Nor(Kx,\Psi)$ with $K\in\cS$, $\Psi\succ0$ and $\vartheta:=(1-\alpha)\tau$ as
\begin{equation}\label{eq:Jent}
J_{\vartheta}(\pi):=C(K,\Psi)-\vartheta\,h\big(\Nor(0,\Psi)\big),
\end{equation}
the average cost minus $\vartheta$ times the entropy of the probing, since every iteration is a generalized majorization--minimization step~\cite{lange2016mm} (Lemma~\ref{lem:mm}).
For $\vartheta>0$, the unique minimizer of $J_\vartheta$ among these policies is $\Nor\big(K^\star x,\frac\vartheta2(S^\star)^{-1}\big)$, so the entropy term changes only its covariance.
This descent gives no rate, and Theorem~\ref{thm:mb} derives one from the gain step of~\eqref{eq:alg}, indexing iterations by $g$.

\begin{theorem}[Known plant]\label{thm:mb}
Let $K_0\in\cS$, $H_0\succ0$, $\alpha\in(0,1]$, and iterate~\eqref{eq:alg} with $S=S_{K_g}$, $E=E_{K_g}$, for any regularization weights $\tau_g>0$, and let $\Psi_g:=\frac{\tau_g}2H_g^{-1}$ for all $g\ge0$.
Let $\bar s:=\|R\|+\|B\|^2C(K_0)/\usig(\Sigma_W)$.
Then every $K_g\in\cS$ and $P_{K_{g+1}}\preceq P_{K_g}-E_{K_g}^\top H_{g+1}^{-1}E_{K_g}$.
\pitem{(a)} If $\alpha<1$, then $C(K_g)-C^\star\le(1-\theta_\alpha)^g(C(K_0)-C^\star)$ with $\theta_\alpha:=\frac{(1-\alpha)\underline{s}}{((1-\alpha)\|H_0\|+\bar s)\gamma^\star}\in(0,1]$, $\Psi_g\succeq\frac{\tau_g}{2}\big(\|H_0\|+\frac{\bar s}{1-\alpha}\big)^{-1}I$, $K_g\to K^\star$ and $H_g\to \tfrac{S^\star}{1-\alpha}$, so that for a constant regularization weight $\Psi_g\to\frac{(1-\alpha)\tau}{2}(S^\star)^{-1}$.
\pitem{(b)} If $\alpha=1$, then $C(K_g)-C^\star\le\big(\frac{c_1}{c_1+g}\big)^{c_2}(C(K_0)-C^\star)$ with $c_1:=1+\|H_0\|/\bar s$, $c_2:=\frac{\underline{s}}{\bar s\gamma^\star}$, and $\frac{\tau_g}{2(\|H_0\|+g\bar s)}I\preceq\Psi_g\preceq\frac{\tau_g}{2(\lambda_{\min}(H_0)+g\lambda_{\min}(R))}I$.
\end{theorem}

\begin{proof}See Appendix~\ref{app:mb}. The proof combines the Loewner-order descent lemma (Lemma~\ref{lem:loewner}) and its consequence for the cost (Lemma~\ref{lem:gn}) with the bound $H_{g+1}\succeq S_{K_g}$.\end{proof}

Theorem~\ref{thm:mb}(b) shows the collapse that tempering prevents. Without tempering, the probing covariance decays like $1/g$ when the regularization weights are bounded above and away from zero.
For $\alpha<1$, the covariance keeps a floor proportional to $\tau_g$, and for a constant regularization weight it tends to $\frac{(1-\alpha)\tau}2(S^\star)^{-1}$, the covariance of the minimizer of $J_\vartheta$, which adds $\tfrac{m(1-\alpha)\tau}{2}$ to the average cost (Proposition~\ref{prop:price}).

\runin{Single steps against whole iterations}
Proposition~\ref{prop:corr} matches one step of KLPI with one step of PGAC, but over many iterations the two methods form their preconditioners and probing covariances differently.
PGAC recomputes its preconditioner $H$ from the current gain alone, and the designer chooses its step size and probing covariance.
KLPI takes $D=\alpha H_t$ from the policy it has just produced, so for $\alpha<1$ its preconditioner remembers past input Hessians with geometrically decaying factors.
The same geometric averaging appears for state-action values in~\cite{vieillard2020leverage} and, read as momentum, in the policy updates of~\cite{mahmoudi2026kl}.
KLPI's probing covariance follows from a chosen regularization weight.
KLPI reproduces PGAC's Gauss--Newton iterates if $\alpha H_t=\frac{1-2\eta}{2\eta}S_{K_t}$ at every $t$, and only if this equality holds at every iteration at which $E_{K_t}$ has full row rank.

\section{The adaptive loop}\label{sec:adaptive}

This section applies KLPI in the adaptive loop, with one iteration per sample in which the plant estimate replaces the plant in the policy evaluation.
Under the data condition and with small updates ($\alpha$ close to one), we prove a guarantee like PGAC's, i.e., the switching closed loop is stable and the cost converges linearly up to a bias, and we show that the probing covariance keeps a floor proportional to the regularization weight.

\begin{algorithm}[t]
\caption{KLPI on the estimated plant}\label{alg:akpi}
\footnotesize
\begin{algorithmic}[1]
\Require offline data with $\Phi_{t_0}\succ0$ and plant estimate $\hat\Theta_{t_0}$; gain $K_{t_0}\in\cS$
\Parameters $H_{t_0}\succ0$; $\alpha\in(0,1]$; regularization weights $\tau_t>0$
\For{$t=t_0,t_0+1,\dots$}
\State $\Psi_t\gets\frac{\tau_t}2H_t^{-1}$; draw $\Omega_t\sim\Nor(0,\Psi_t)$ \Comment{probing}
\State apply $U_t=K_tX_t+\Omega_t$ and observe $X_{t+1}$
\State update $\hat\Theta_{t+1}$ by recursive least squares
\State $\hat S_t,\hat E_t\gets$ input Hessian and slope of $K_t$ on $(\hat A_{t+1},\hat B_{t+1})$
\State $H_{t+1}\gets\alpha H_t+\hat S_t$ \Comment{preconditioner}
\State $K_{t+1}\gets K_t-H_{t+1}^{-1}\hat E_t$ \Comment{gain update}
\EndFor
\end{algorithmic}
\end{algorithm}

Algorithm~\ref{alg:akpi} applies~\eqref{eq:alg} once per sample and needs the same computation per sample as PGAC's Gauss--Newton step.
It differs from that step in two ways.
Its preconditioner $H_{t+1}$ accumulates the estimated input Hessians of past samples instead of using the current one alone, and the tempering exponent $\alpha$ of Section~\ref{sec:alg} acts in this sum as a forgetting factor.
Its probing covariance $\frac{\tau_t}2H_t^{-1}$ follows from this preconditioner and the regularization weight, instead of being chosen separately.

\begin{theorem}[Adaptive loop]\label{thm:adaptive}
Let Assumption~\ref{ass:pe} hold, $K_{t_0}\in\cS$ and $\alpha\in(0,1)$.
There are constants $\nu_1,\dots,\nu_5>0$ with $\nu_3\ge1\ge\nu_4$ and $\sa\ge\|R\|+\|B\|^2(C(K_{t_0})+2)/\usig(\Sigma_W)$, depending only on $(A,B,Q,R,\Sigma_W,C(K_{t_0}))$, such that, if
\pitem{(C1)} $1-\alpha\le\nu_2$,
\pitem{(C2)} $\SNR_t\ge\nu_1$ for all $t\ge t_0$, and
\pitem{(C3)} $\frac{\underline{s}}{2(1-\alpha)}I\preceq H_{t_0}\preceq\frac{\sa}{1-\alpha}I$,
\par\noindent then, for all regularization weights $\tau_t>0$, the estimated input Hessians and slopes $\hat S_t,\hat E_t$ of Algorithm~\ref{alg:akpi} are well defined at every sample, the gains $K_t$ are sequentially stable,
\begin{equation}\label{eq:statebound}
\|X_t\|\le \nu_3\big(1-\tfrac{\nu_4}{2}\big)^{t-t_0}\|X_{t_0}\|+\tfrac{2\nu_3}{\nu_4}\max_{t_0\le i<t}\|B\Omega_i+W_i\|,
\end{equation}
and, with $q_\alpha:=1-\frac{(1-\alpha)\underline{s}}{2\sa\gamma^\star}$,
\begin{equation}\label{eq:adrate}
\begin{aligned}
C(K_t)-C^\star\le{}&q_\alpha^{t-t_0}\big(C(K_{t_0})-C^\star\big)\\
&+(1-\alpha)\nu_5\textstyle\sum_{i=t_0+1}^{t}q_\alpha^{t-i}\,\SNR_{i}^{-1}.
\end{aligned}
\end{equation}
Under (C2), (C3) holds for $H_{t_0}=\tfrac{\hat S_{\mathrm{off}}}{1-\alpha}$, where $\hat S_{\mathrm{off}}$ is the input Hessian of $K_{t_0}$ on $\hat\Theta_{t_0}$.
\end{theorem}

\begin{proof}See Appendix~\ref{app:adaptive}. The proof follows the proofs of~\cite[Theorems~1 and~3]{zhao2025pgac}. It compares the computed update with the exact update that uses the true slope, and bounds their difference by a multiple of $1/\SNR_{t+1}$, where $\SNR_{t+1}$ measures the data behind the plant estimate $\hat\Theta_{t+1}$. Two ingredients are new: a Loewner-order descent lemma for matrix preconditioners (Lemma~\ref{lem:loewner}), and an induction that keeps $H_t$ in the band of (C3).\end{proof}
Condition (C1) is the counterpart of PGAC's small step size, and the probing enters only through (C2) and $B\Omega_i$ in~\eqref{eq:statebound}.
Appendix~\ref{app:adaptive} gives every constant of Theorem~\ref{thm:adaptive} in closed form, except a supremum of the slope norm $\|E_K\|$ inside $\nu_2$.

\begin{corollary}[Probing in the adaptive loop]\label{cor:floor}
Under the conditions of Theorem~\ref{thm:adaptive}, for all $t\ge t_0$, the probing covariance $\Psi_t=\frac{\tau_t}2H_t^{-1}$ of Algorithm~\ref{alg:akpi} satisfies
\begin{equation}\label{eq:floor}
\tfrac{(1-\alpha)\tau_t}{2\sa}I\preceq\Psi_t\preceq\tfrac{(1-\alpha)\tau_t}{\underline{s}}I.
\end{equation}
\end{corollary}
\begin{proof}Starting from (C3), the induction of Appendix~\ref{app:adaptive} applies Lemma~\ref{lem:onestep}(i) at every $t\ge t_0$ and so keeps $H_t$ in the band of (C3).
Inverting this band gives~\eqref{eq:floor}.\end{proof}
Hence, as on a known plant (Theorem~\ref{thm:mb}(a)), the probing does not die out in the adaptive loop when the regularization weights are bounded away from zero.
We leave open whether this floor on the probing covariance implies the data condition (C2).

\section{Comparing KLPI and PGAC in Simulation}\label{sec:cmp}\label{sec:sim}

This section compares in simulation the two design choices in which KLPI differs from PGAC's Gauss--Newton step.
The first is the probing \textbf{\emph{shape}}~\eqref{eq:amountshape}, since KLPI probes with a covariance proportional to the inverse of its preconditioner.
The second is the \textbf{\emph{averaging}} of past estimated input Hessians in this preconditioner.
We define the compared quantity (Section~\ref{sec:cmp-cost}), describe the simulation setup (Section~\ref{sec:sim-setup}), and report the effect of each choice (Sections~\ref{sec:sim-shape} and~\ref{sec:sim-avg}).

\subsection{The Compared Quantity}\label{sec:cmp-cost}

This subsection defines the ratio we compare, from a split of the regret, and the input directions along which probing shapes differ.

\runin{Regret and its parts}
We measure each method by its regret over $T$ samples,
\begin{equation}\label{eq:regdef}
\mathrm{Reg}_T:=\EE_{\PP}\Big[\textstyle\sum_{t<T}c(X_t,U_t)\Big]-TC^\star,
\end{equation}
where $\PP$ is the law of the closed loop under the method.
For a method that applies $U_t=K_tX_t+\Omega_t$ with $\Omega_t\sim\Nor(0,\Psi_t)$, we define two costs at each sample $t$: the probing cost $\mathcal P_t$, and the gain-error cost $\mathcal E_t$, which is the cost caused by the suboptimal gain $K_t$:
\begin{equation}\label{eq:pecost}
\mathcal P_t:=\Tr(S^\star\Psi_t),\ \
\mathcal E_t:=X_t^\top(K_t-K^\star)^\top S^\star(K_t-K^\star)X_t.
\end{equation}
Proposition~\ref{prop:regret} splits the regret into a boundary term and the sums of these two costs.

\begin{proposition}[Cost split]\label{prop:regret}
Let $U_t=K_tX_t+\Omega_t$, where the gain $K_t$ and the probing covariance $\Psi_t$ are functions of the data $\mathcal I_t:=(X_0,U_0,\dots,U_{t-1},X_t)$ observed up to time $t$.
Let the probing perturbation $\Omega_t$ have the conditional distribution $\Nor(0,\Psi_t)$ given $\mathcal I_t$, and let the process noise $W_t$ be independent of $(\mathcal I_t,\Omega_t)$.
If, under the closed-loop law $\PP$, $\EE_{\PP}\|X_t\|^2$, $\EE_{\PP}\|K_tX_t\|^2$ and $\EE_{\PP}\Tr\Psi_t$ are finite for $t\le T$, then
\begin{equation}\label{eq:regret}
\resizebox{0.91\columnwidth}{!}{$\begin{aligned}
\EE_{\PP}\Big[\textstyle\sum_{t<T}c(X_t,U_t)\Big]={}&TC^\star+\EE_{\PP}\big[X_0^\top P^\star X_0-X_T^\top P^\star X_T\big]\\
&+\textstyle\sum_{t<T}\EE_{\PP}[\mathcal P_t]+\sum_{t<T}\EE_{\PP}[\mathcal E_t].
\end{aligned}$}
\end{equation}
\end{proposition}

\begin{proof}See Appendix~\ref{app:adaptive}.\end{proof}

The probing cost $\mathcal P_t$ does not depend on the applied gain, so two methods that probe with the same covariance pay the same probing cost.

\runin{Equal probing cost}
Probing makes the data more informative, so a method that probed more could improve its gains for that reason alone.
To make the two methods comparable, we therefore give both the same probing cost at every sample.
We fix this common cost through a probing budget $c$ per sample, chosen by the designer (Section~\ref{sec:sim-setup}).
The methods start after $t_0$ offline samples (Section~\ref{sec:sim-setup}), so we apply Proposition~\ref{prop:regret} from the first online sample $t_0$ on.
The summed probing costs of the two methods in~\eqref{eq:regret} are then equal, and so are the first parts $\EE_{\PP}[X_{t_0}^\top P^\star X_{t_0}]$ of their boundary terms, since both methods start from the same state $X_{t_0}$.
The last part $\EE_{\PP}[X_T^\top P^\star X_T]$ depends only on the final state and does not grow with $T$ when the gains stabilize the plant.
The summed gain-error cost, in contrast, adds a term at every sample, so for long runs the regrets differ mainly through the summed gain-error costs.
We run each method with many random seeds $s$, and we compare two methods $\mathrm a$ and $\mathrm b$ by the ratio of their summed gain-error costs:
\begin{equation}\label{eq:sim-ratio}
r:=\textstyle\sum_s\sum_{t=t_0}^{t_0+T-1}\mathcal E^{\mathrm a}_{t,s}\big/\sum_s\sum_{t=t_0}^{t_0+T-1}\mathcal E^{\mathrm b}_{t,s}.
\end{equation}
We call a difference clear when the $95\%$ paired-bootstrap interval of $r$ over the seeds excludes one.
Since the designer of each method chooses how much to probe, we repeat the comparison at the budgets $c/C^\star\in\{1/64,1/16,1/4,1/2,1,2\}$.
We do not consider larger budgets, at which the input would serve mainly to identify the plant rather than to control it.

\runin{Cheap and expensive input directions}
We call the unit eigenvectors of $S^\star$ the \textbf{\emph{input directions}}.
By~\eqref{eq:pecost}, probing with variance $\psi$ along an input direction with eigenvalue $\lambda$ costs $\lambda\psi$, so different probing shapes can divide the probing cost differently among the input directions.
We call an input direction cheap when its eigenvalue lies below the mean eigenvalue of $S^\star$, and expensive otherwise.
We measure how much the input directions differ in cost by the quantity $\varrho:=\lambda_{\max}(S^\star)/\lambda_{\min}(S^\star)$, which we call the \textbf{\emph{input spread}}.

\begin{remark}[KLPI's probing shape against isotropic probing]\label{rem:shape}
With a known plant, $\alpha<1$ and a constant regularization weight, KLPI's probing covariance approaches a multiple of $(S^\star)^{-1}$ (Theorem~\ref{thm:mb}(a)).
This covariance gives an input direction with eigenvalue $\lambda$, a variance proportional to $1/\lambda$, and hence an equal share of the probing cost.
Isotropic probing gives every input direction the same variance, and hence a share proportional to $\lambda$.
At the same probing cost, the shape $(S^\star)^{-1}$ therefore probes every cheap input direction more than isotropic probing does, and every expensive one less.
The two shapes coincide when $\varrho=1$.
\end{remark}

\subsection{Simulation Setup}\label{sec:sim-setup}

This subsection describes the four variants of KLPI and PGAC, how much they probe, the plants and the two ways a run starts.
The code of all simulations is available at~\repolink.

{%
\runin{Compared variants}
Both methods update a recursive least-squares estimate of the plant at every sample and compute their gain update on this plant estimate.
KLPI's gain update is a Gauss--Newton step with step size $\eta=\tfrac{1-\alpha}{2}$ and an averaged input Hessian (Section~\ref{sec:alg}).
We let PGAC take Gauss--Newton steps with the same step size, and we set $\alpha=0.9$.
To separate the two design choices, each method probes either isotropically, as in~\cite{zhao2025pgac}, or with a covariance proportional to the inverse of its own preconditioner.
The four combinations are the variants of Table~\ref{tab:variants}.
\runin{Amount of probing}
For KLPI, we use the constant regularization weight $\tau=2c/(m(1-\alpha))$, for which Theorem~\ref{thm:mb}(a) gives a probing cost $\mathcal P_t$ that tends to the budget $c$.
The other variants probe with $\Psi_t=\sigma_t^2\Lambda_t$, where $\Lambda_t$ is the shape in Table~\ref{tab:variants}.
We choose $\sigma_t^2$ so that each of these variants pays the same probing cost~\eqref{eq:pecost} as KLPI at every sample, where KLPI probes with $\Psi^{\mathrm{KLPI}}_t=\tfrac\tau2(H^{\mathrm{KLPI}}_t)^{-1}$:
\begin{equation}\label{eq:sim-arms}
\sigma_t^2\Tr(S^\star\Lambda_t)=\Tr\big(S^\star\Psi^{\mathrm{KLPI}}_t\big).
\end{equation}
This matching is the only computation that uses the true input Hessian $S^\star$, and no gain update uses it.

\begin{table}[t]
\caption{Gain update and probing of the four variants}
\label{tab:variants}
\centering
\scriptsize
\renewcommand{\arraystretch}{1.15}
\begin{tabular}{@{}lcc@{}}
\toprule
Variant & Preconditioner $H_{t+1}$ & Probing shape $\Lambda_t$\\
\midrule
PGAC-iso & $\hat S_t/(1-\alpha)$ & $I$\\
PGAC-shaped & $\hat S_t/(1-\alpha)$ & $H_t^{-1}$\\
KLPI-iso & $\alpha H_t+\hat S_t$ & $I$\\
KLPI & $\alpha H_t+\hat S_t$ & $H_t^{-1}$\\
\bottomrule
\end{tabular}

\vspace{2pt}
\parbox{\columnwidth}{\scriptsize Every variant starts from $H_{t_0}=\hat S_{\mathrm{off}}/(1-\alpha)$ (Theorem~\ref{thm:adaptive}) and updates its gain by $K_{t+1}=K_t-H_{t+1}^{-1}\hat E_t$.}
\end{table}
\begin{table}[t]
\caption{Plants used in the simulations}
\label{tab:sim-plants}
\centering
\scriptsize
\setlength{\tabcolsep}{1.6pt}
\renewcommand{\arraystretch}{1.15}
\begin{tabular}{@{}lcllll@{}}
\toprule
Plant & $n$ & $A$ & $B$ & $R$ & $\varrho$\\
\midrule
P1 & $4$ & \multicolumn{2}{l}{benchmark of DeePO~\cite[Sec.~VI]{zhao2025deepo}} & $I_2$ & $4.2$\\
P2 & $2$ & $\left[\begin{smallmatrix}0.9&0.3\\0&0.95\end{smallmatrix}\right]$ & $I_2$ & $\operatorname{diag}(1,100)$ & $43$\\
P3 & $2$ & as P2 & $\operatorname{diag}(1,9.927)$ & $I_2$ & $43$\\
P4 & $2$ & $\operatorname{diag}(0,0.9)$ & $\operatorname{diag}(2,1)$ & $\operatorname{diag}(100,1)$ & $42$\\
random (40) & $3$ & Gaussian$^{\mathrm a}$ & Gaussian & $\operatorname{diag}(1,\varphi)$ & varies\\
\bottomrule
\end{tabular}

\vspace{2pt}
\parbox{\columnwidth}{\scriptsize All plants have two inputs, $Q=I_n$ and $\Sigma_W=I_n$.
$^{\mathrm a}$Rescaled to a spectral radius drawn uniformly from $[0.5,1.05]$.
For the random plants, twenty have $\varphi=1$ and twenty have $\log_{10}\varphi$ drawn uniformly from $[0,2]$.}
\end{table}

\begin{figure*}[t]
\centering
\includegraphics[width=\textwidth]{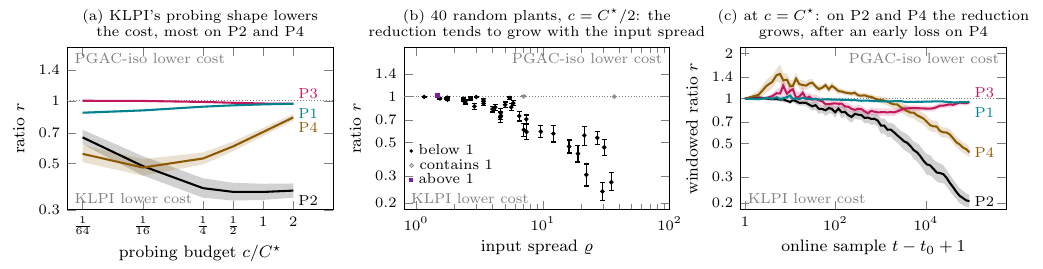}
\setlength{\abovecaptionskip}{0pt}
\vspace{-16pt}
\caption{Ratio $r$ in~\eqref{eq:sim-ratio} of KLPI to PGAC-iso from the standard start, with $95\%$ bootstrap intervals over the seeds ($480$ per plant for P1--P4, $50$ per random plant): (a) whole runs on P1--P4; (b) whole runs on the random plants, $c=C^\star/2$; (c) windowed ratio on P1--P4, $c=C^\star$, with windows evenly spaced on the logarithmic axis.}
\label{fig:sim-shape}
\end{figure*}

\runin{Plants}
We test the variants on four plants, each chosen to bring out one effect, and on forty random plants that guard against conclusions specific to those four (Table~\ref{tab:sim-plants}).
P1 is the DeePO benchmark~\cite{zhao2025deepo}, whose input spread is the smallest of P1--P4, so it serves as a control case.
P2 and P3 have the same $A$ and a large input spread, obtained from the input weight $R$ on P2 and from the input matrix $B$ on P3.
P2 tests whether KLPI's probing shape helps when the input directions differ strongly in cost, and P3 tests whether the source of the input spread matters.
P4 has nearly the input spread of P2, but its $A$, $B$ and $R$ are diagonal, so each input drives its own state and the input directions are the two inputs.
Its cheap input drives a state that decays slowly, whereas its expensive input affects its state for one sample only.
Compared with isotropic probing, KLPI's probing shape moves probing from the expensive input to the cheap one (Remark~\ref{rem:shape}), so on P4 its probing enlarges the slowly decaying state.
Since the gain-error cost~\eqref{eq:pecost} is quadratic in the state, P4 tests whether KLPI's probing shape still helps when its probing makes a gain error costlier.

\runin{Starts}
Every run begins with $t_0=100$ offline samples under independent standard Gaussian inputs, to which we fit a least-squares estimate of the plant, followed by $T=10^5$ online samples.
In the \textbf{\emph{standard start}}, the variants begin, as in~\cite{zhao2025pgac}, with the certainty-equivalent gain of this plant estimate.
The \textbf{\emph{poor-gain start}}, used in Section~\ref{sec:sim-avg}, differs only in the initial gain, which is the certainty-equivalent gain of a plant estimate fitted to the first eight offline samples only.
The gain updates still use the plant estimate fitted to all $t_0$ offline samples, so the poor initial gain is paired with a much more accurate plant estimate.
Although this pairing is not a practical scenario, we use it as a controlled test of the averaging when the gain has much to improve.
\par}

\subsection{The Probing Shape}\label{sec:sim-shape}

We compare KLPI with PGAC-iso from the standard start by the ratio $r$ in~\eqref{eq:sim-ratio}, which Fig.~\ref{fig:sim-shape} presents in three views.
In Fig.~\ref{fig:sim-shape}(a), $r$ covers whole runs on P1--P4 at every budget.
In Fig.~\ref{fig:sim-shape}(b), $r$ covers whole runs on the forty random plants at $c=C^\star/2$ and is plotted against the input spread.
Fig.~\ref{fig:sim-shape}(c) follows P1--P4 over time at $c=C^\star$, by the windowed ratio, which is $r$ computed over a short window of online samples.
KLPI and PGAC-iso differ in both design choices.
However, KLPI and PGAC-shaped, which differ only in the averaging, have summed gain-error costs within $0.6\%$ of each other on P1--P4 at every budget and on every random plant at $c=C^\star/2$.
Fig.~\ref{fig:sim-shape} therefore omits PGAC-shaped, whose ratios in panels (a) and (b) would coincide with those of KLPI.
The differences in summed gain-error cost between KLPI and PGAC-iso in Fig.~\ref{fig:sim-shape} thus come from KLPI's probing shape, and PGAC-shaped shows the same differences against PGAC-iso.

{%
KLPI's probing shape lowers the summed gain-error cost on nearly every tested plant and budget, most on P2 and P4 (Fig.~\ref{fig:sim-shape}(a)).
On the random plants, $r$ tends to be smaller on plants with a larger input spread (Fig.~\ref{fig:sim-shape}(b)).
Intuitively, the gain update multiplies the estimated slope by the inverse preconditioner, which is largest along the cheap input directions, so errors in the plant estimate move the gain update most along these directions.
KLPI's probing shape excites these directions more than isotropic probing does (Remark~\ref{rem:shape}), and thus makes the plant estimate more accurate along these directions.
The windowed ratio on P2 falls nearly steadily over the run (Fig.~\ref{fig:sim-shape}(c)).

The probing shape helps less on P3, and costs more early in the run on P3 and P4, and briefly on P1.
P3 has the same $A$ and the same input spread as P2, but its spread comes from $B$ instead of $R$.
On P3 the ratio $r$ over whole runs stays close to one at every budget (Fig.~\ref{fig:sim-shape}(a)).
On P3 the windowed ratio exceeds one for the first few dozen samples, falls to about $0.8$ and then rises back towards one, whereas on P2 it keeps falling (Fig.~\ref{fig:sim-shape}(c)).
The input spread alone therefore does not determine how much the probing shape helps.
On P4, KLPI pays more than PGAC-iso early in the run and less afterwards (Fig.~\ref{fig:sim-shape}(c)).
We attribute this early loss to KLPI's stronger probing of the cheap input of P4, which enlarges the slowly decaying state, and with it the gain-error cost~\eqref{eq:pecost}, before the gain has improved.
In runs at $c=2C^\star$, which Fig.~\ref{fig:sim-shape}(c) does not show, the early loss is larger, which is consistent with the rise of the ratio on P4 in Fig.~\ref{fig:sim-shape}(a).
Over whole runs, KLPI's probing shape clearly raises the summed gain-error cost only on one random plant, and only slightly (Fig.~\ref{fig:sim-shape}(b)).

Overall, KLPI's probing shape lowers the summed gain-error cost over whole runs on nearly every tested plant and budget, but by amounts that differ widely between plants, and not at every stage of a run.
Intuitively, the probing covariance acts on the summed gain-error cost in two ways.
It sets how strongly each input direction is excited, and thereby along which input directions the plant estimate becomes accurate, which shapes the gain updates.
The probing also drives the state through $B$, and the gain-error cost~\eqref{eq:pecost} grows quadratically with the state.
How the two effects combine depends on the plant, the budget and the stage of the run, and they can pull in opposite directions, as the early loss on P4 shows (Fig.~\ref{fig:sim-shape}(c)).
A model that predicts when a probing shape helps and when it hurts would therefore help to design the probing.
We do not develop such a model in this paper and leave it to future work.
\par}

\begin{figure}[t]
\centering
\includegraphics[width=\columnwidth]{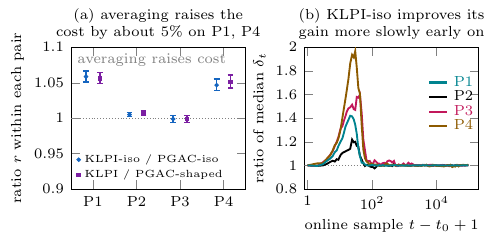}
\setlength{\abovecaptionskip}{0pt}
\vspace{-16pt}
\caption{Effect of the averaging from the poor-gain start at $c=C^\star$ ($480$ seeds per plant): (a) the ratio $r$ in~\eqref{eq:sim-ratio} within each pair that differs only in the averaging, with $95\%$ bootstrap intervals; (b) the ratio of the median relative gap $\delta_t$ of KLPI-iso to that of PGAC-iso, averaged over short windows of online samples.}
\label{fig:sim-avg}
\end{figure}

\subsection{The Averaging}\label{sec:sim-far}\label{sec:sim-avg}

This subsection asks how the averaging affects the gain-error cost from each of the two starts.
We use the two pairs of variants of Table~\ref{tab:variants} that differ only in the averaging: KLPI-iso and PGAC-iso, and KLPI and PGAC-shaped.
To see how fast each method improves its gain, we follow the relative gap of the gain, $\delta_t:=(C(K_t)-C^\star)/C^\star$, where $C(K_t)$ is the average cost of the gain $K_t$ without probing (Section~\ref{sec:setting}).

{%
From the standard start, the ratio within each of the two pairs is close to one on P1--P4 at every budget, so Fig.~\ref{fig:sim-avg} shows only the poor-gain start.
Intuitively, the averaged and the current input Hessians differ much only while the input Hessian changes quickly, as it does while a poor gain improves.

From the poor-gain start, the averaging raises the gain-error cost slightly on some plants, and the probing shape barely alters this effect (Fig.~\ref{fig:sim-avg}(a)).
Early in the run, KLPI-iso improves the gain more slowly than PGAC-iso, and afterwards the two relative gaps agree (Fig.~\ref{fig:sim-avg}(b)).
Intuitively, as a poor gain improves, its input Hessian $S_K$ shrinks with $P_K$ (Lemma~\ref{lem:loewner}), so an average over past input Hessians stays above the current one.
The averaged preconditioner is then larger than PGAC's, which shortens the gain step.

\par}

\section{Conclusion}\label{sec:concl}

We approached adaptive LQR through a KL-regularized formulation, in which the probing is the randomness of a Gaussian policy and its cost enters the objective.
This formulation offers one way to choose the probing covariance, which earlier works on policy-gradient adaptive control~\cite{zhao2025deepo,zhao2025pgac} leave as a design parameter.
Each policy improvement yields a preconditioned gain update together with a probing covariance proportional to the inverse preconditioner (Proposition~\ref{prop:soft}).
With a suitable reference centred at the current gain, this gain update is PGAC's Gauss--Newton or natural-gradient step for any step size below a bound (Proposition~\ref{prop:corr}).
We repeat the policy improvement against a tempered copy of the current policy, which gives KLPI, and show that with a known plant and $\alpha<1$ the cost of its gain converges linearly and its probing covariance keeps a floor proportional to the regularization weight (Theorem~\ref{thm:mb}).
In the adaptive loop, with small updates, a suitably initialized preconditioner and the data condition, KLPI keeps the closed loop stable and its cost converges linearly up to a bias (Theorem~\ref{thm:adaptive}), and its probing covariance keeps a floor proportional to the regularization weight (Corollary~\ref{cor:floor}).
In simulations at equal probing cost, KLPI's probing shape, which PGAC can adopt, lowered the summed gain-error cost compared with isotropic probing on nearly all tested plants and budgets, by amounts that varied widely between plants.
The probing shape raised the gain-error cost early in the run on some plants, and slightly over whole runs on one random plant.
KLPI's averaging of past input Hessians changed the gain-error cost little, except in a controlled test that starts from a poor gain, where it slowed the early improvement of the gain.
We leave to future work a way to predict when a probing shape helps and when it hurts, and a regret bound, including which schedules of the regularization weight and of $\alpha$ make the regret sublinear.

\section*{Acknowledgment}
The authors used Claude (Anthropic), a large language model, in preparing this work.
It drafted the proofs, the simulation code and parts of the text under the direction of the authors, who checked the results and take full responsibility for the content of this paper.

\bibliographystyle{IEEEtran}
\bibliography{references}

\onecolumn
\renewcommand{\baselinestretch}{1.0}\selectfont
\begin{adjustwidth}{0.6in}{0.6in}
\useRomanappendicesfalse
\appendices

\section*{Overview of the Appendices}

We collect the proofs in four appendices, and Fig.~\ref{fig:appmap} shows which results each proof uses.

Appendix~\ref{app:prelim} collects the facts about the LQR that the later proofs use, and its first four lemmas build on one another.
The cost-difference identity (Lemma~\ref{lem:costdiff}) expresses the change of the cost between two gains, and gradient dominance (Lemma~\ref{lem:gd}) follows from this identity.
The expansion~\eqref{eq:defect} in Lemma~\ref{lem:costdiff} also gives a descent lemma in the Loewner order (Lemma~\ref{lem:loewner}).
This lemma shows that a gain step with a large enough preconditioner keeps the closed loop stable and decreases the solution $P_K$ of its Lyapunov equation.
Combined with gradient dominance, Lemma~\ref{lem:loewner} shows that each of these steps shrinks the gap $C(K)-C^\star$ to the optimal cost by a guaranteed fraction (Lemma~\ref{lem:gn}).
Lemma~\ref{lem:bounds} bounds $P_K$, $\Sigma_K$ and $K$ for all gains of bounded cost, which gives the constants of Appendices~\ref{app:mb} and~\ref{app:adaptive}.
Lemma~\ref{lem:lyap} bounds how the solution of a Lyapunov equation changes with its data, which Appendix~\ref{app:adaptive} uses to compare the plant estimate with the true plant.

Appendix~\ref{app:kl} proves Propositions~\ref{prop:soft} to~\ref{prop:corr} on one KL-regularized policy improvement, and the proof of Proposition~\ref{prop:soft} rests on the Gibbs variational identity (Lemma~\ref{lem:gibbs}).
Lemma~\ref{lem:mm} then shows that the iterations of KLPI do not increase the fixed objective of Section~\ref{sec:alg}.

Appendix~\ref{app:mb} proves Theorem~\ref{thm:mb} for a known plant by applying Lemmas~\ref{lem:loewner} and~\ref{lem:gn} at every iteration, with the bounds of Lemma~\ref{lem:bounds}.

Appendix~\ref{app:adaptive} proves Theorem~\ref{thm:adaptive} on the adaptive loop and Proposition~\ref{prop:regret} on the cost split.
Fact~\ref{f:pert} bounds the error of the objects of a gain computed on the plant estimate, and Fact~\ref{f:lip} bounds how these objects change with the gain.
Facts~\ref{f:strong} and~\ref{f:state} give the stability of the closed loop and a bound on its state.
With these facts, Lemma~\ref{lem:hessband} keeps the estimated input Hessian in a fixed band, and Lemma~\ref{lem:onestep} bounds the cost after one step of Algorithm~\ref{alg:akpi}.
An induction over the samples combines these results into Theorem~\ref{thm:adaptive}.

\begin{figure}[h]
\centering
\includegraphics[width=\linewidth]{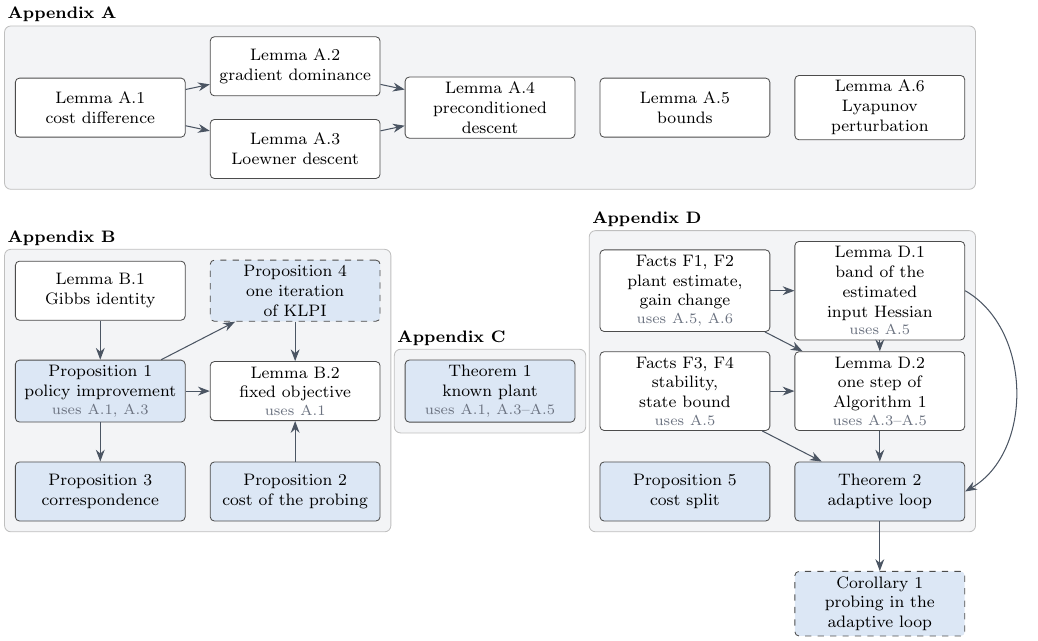}
\caption{Dependencies among the results of the appendices. An arrow from one result to another means that the proof of the second uses the first, and the grey line in a box lists the lemmas of Appendix~\ref{app:prelim} used in the proof of that result. Shaded boxes are results stated in the body of the paper, and dashed boxes are proved there.}
\label{fig:appmap}
\end{figure}

\section{Preliminaries}\label{app:prelim}

This appendix collects the facts about the LQR that the proofs use: a cost-difference identity, gradient dominance, a descent lemma for preconditioned steps in the Loewner order with its consequence for the cost, bounds on gains of bounded cost, and a perturbation bound for Lyapunov equations.
The first four lemmas build on one another, and the last two supply bounds for the later appendices.
For $v\ge C^\star$ let $\cS(v):=\{K\in\cS:C(K)\le v\}$.
Lemmas~\ref{lem:costdiff} and~\ref{lem:gd} are proved in~\cite{fazel2018global} for a noise-free plant with a random initial state.
Their statements hold in the setting of this paper with $\Sigma_W$ in place of the initial-state covariance of~\cite{fazel2018global}, because in both settings the cost is the trace of $P_K$ times this covariance, and $P_K$ and $\Sigma_K$ solve Lyapunov equations of the same form.

\begin{applemma}[{Cost difference~\cite[Lemmas~10 and~12]{fazel2018global}}]\label{lem:costdiff}
For $K,K'\in\cS$ and $\Delta:=K'-K$,
\begin{equation}\label{eq:costdiff}
C(K')-C(K)=\Tr\big(\Sigma_{K'}(\Delta^\top E_K+E_K^\top\Delta+\Delta^\top S_K\Delta)\big),
\end{equation}
where, with $P:=P_K$, the matrix in the trace is
\begin{equation}\label{eq:defect}
\begin{aligned}
Y&:=Q+K'^\top RK'+A_{K'}^\top PA_{K'}-P\\
&\phantom{:}=\Delta^\top E_K+E_K^\top\Delta+\Delta^\top S_K\Delta.
\end{aligned}
\end{equation}
The expansion~\eqref{eq:defect} holds for every $K'\in\mathbb R^{m\times n}$, since it involves $P_K$ but neither $P_{K'}$ nor $\Sigma_{K'}$.
In particular,~\eqref{eq:costdiff} with $K^\star$ in place of $K$ and $K$ in place of $K'$ gives, since $E_{K^\star}=0$, $C(K)-C^\star=\Tr(\Sigma_K(K-K^\star)^\top S^\star(K-K^\star))$, which vanishes only at $K=K^\star$.
Moreover, for $K\in\cS$, $P_K-P^\star=A_K^\top(P_K-P^\star)A_K+(K-K^\star)^\top S^\star(K-K^\star)$, so $P_K\succeq P^\star$.
\end{applemma}

The identity for $P_K-P^\star$ in Lemma~\ref{lem:costdiff} is not stated in~\cite{fazel2018global}; it follows from~\eqref{eq:defect}.
Taking $K^\star$ in place of $K$, $K$ in place of $K'$ and $E_{K^\star}=0$ in~\eqref{eq:defect} gives $Q+K^\top RK+A_K^\top P^\star A_K-P^\star=(K-K^\star)^\top S^\star(K-K^\star)$, and subtracting this identity from the Lyapunov equation $P_K=Q+K^\top RK+A_K^\top P_KA_K$ gives the identity for $P_K-P^\star$.
This identity is a Lyapunov equation with the Schur matrix $A_K$ and a positive semidefinite constant term, so $P_K-P^\star\succeq0$.

In~\cite{fazel2018global}, gradient dominance is derived from Lemma~\ref{lem:costdiff}.
Completing the square in the gain change $\Delta$ bounds the matrix $Y$ of~\eqref{eq:defect} from below by $-E_K^\top S_K^{-1}E_K$, and~\eqref{eq:costdiff} with $K^\star$ in place of $K'$ then gives the bound.

\begin{applemma}[{Gradient dominance~\cite[Lemma~11]{fazel2018global}}]\label{lem:gd}
For $K\in\cS$, $C(K)-C^\star\le\|\Sigma^\star\|\Tr\big(E_K^\top S_K^{-1}E_K\big)$.
\end{applemma}

\begin{applemma}[Loewner descent]\label{lem:loewner}
Let $K\in\cS$, $S:=S_K$, $E:=E_K$, $\chi\in(0,2)$, $M\succ0$ with $MSM\preceq\chi M$, and $K':=K-ME$.
Then $K'\in\cS$ and $P_{K'}\preceq P_K-(2-\chi)E^\top ME$.
\end{applemma}

\begin{proof}
Let $F(P'):=Q+K'^\top RK'+A_{K'}^\top P'A_{K'}$ be the Lyapunov map of $K'$.
We first evaluate $F(P_K)-P_K$, and then deduce that $K'$ is stabilizing and that $P_{K'}$ lies below $P_K$.
With $P=P_K$, the matrix $Y$ in~\eqref{eq:defect} is $F(P_K)-P_K$, and the gain change is $\Delta=-ME$.
Since $M$ is symmetric, inserting $\Delta=-ME$ into the right side of~\eqref{eq:defect} gives $-2E^\top ME+E^\top MSME=-E^\top(2M-MSM)E$.
Since $MSM\preceq\chi M$, we have $2M-MSM\succeq(2-\chi)M$, and therefore
\begin{equation}\label{eq:loew1}
F(P_K)-P_K=-E^\top(2M-MSM)E\preceq-(2-\chi)E^\top ME.
\end{equation}
Next we show that $K'$ is stabilizing.
Since $F(P_K)\preceq P_K$ by~\eqref{eq:loew1}, we have $A_{K'}^\top P_KA_{K'}-P_K\preceq-Q\prec0$ with $P_K\succ0$.
For an eigenvector $v$ of $A_{K'}$ with eigenvalue $\zeta_v$, this matrix inequality gives $(|\zeta_v|^2-1)v^\ast P_Kv<0$, so $|\zeta_v|<1$ and $K'\in\cS$.
It remains to compare $P_{K'}$ with $P_K$.
Since $P_{K'}=F(P_{K'})$ and $F(P_K)-F(P_{K'})=A_{K'}^\top(P_K-P_{K'})A_{K'}$, the difference $P_K-P_{K'}$ satisfies
\begin{equation*}
\begin{aligned}
P_K-P_{K'}&=\big(P_K-F(P_K)\big)+\big(F(P_K)-F(P_{K'})\big)\\
&=\big(P_K-F(P_K)\big)+A_{K'}^\top(P_K-P_{K'})A_{K'},
\end{aligned}
\end{equation*}
a Lyapunov equation with the Schur matrix $A_{K'}$.
Its solution is $\sum_{k\ge0}(A_{K'}^\top)^k(P_K-F(P_K))A_{K'}^k$, whose terms are positive semidefinite by~\eqref{eq:loew1}, so the solution is at least its term $k=0$, which is $P_K-F(P_K)$.
Therefore~\eqref{eq:loew1} gives $P_K-P_{K'}\succeq(2-\chi)E^\top ME$.
\end{proof}

\begin{applemma}[Preconditioned descent]\label{lem:gn}
Under the hypotheses of Lemma~\ref{lem:loewner}, with $\varkappa:=\lambda_{\min}(S^{1/2}MS^{1/2})$ and $\gamma^\star$ as in Section~\ref{sec:setting},
\begin{equation}\label{eq:gnrate}
C(K')-C^\star\le\big(1-(2-\chi)\varkappa/\gamma^\star\big)\big(C(K)-C^\star\big),
\end{equation}
and $0<(2-\chi)\varkappa\le1\le\gamma^\star$, so the factor lies in $[0,1)$.
\end{applemma}

\begin{proof}
Multiplying the Loewner decrease of Lemma~\ref{lem:loewner} by $\Sigma_W$ and taking traces gives
\begin{equation*}
\begin{aligned}
C(K)-C(K')&=\Tr\big((P_K-P_{K'})\Sigma_W\big)\\
&\ge(2-\chi)\Tr(\Sigma_WE^\top ME)\\
&\ge(2-\chi)\usig(\Sigma_W)\Tr(E^\top ME)\\
&\ge(2-\chi)\varkappa\,\usig(\Sigma_W)\Tr(E^\top S^{-1}E),
\end{aligned}
\end{equation*}
where the last inequality uses $M\succeq\varkappa S^{-1}$, which follows from $S^{1/2}MS^{1/2}\succeq\varkappa I$.
Lemma~\ref{lem:gd} together with $\gamma^\star=\|\Sigma^\star\|/\usig(\Sigma_W)$ gives $\usig(\Sigma_W)\Tr(E^\top S^{-1}E)\ge(C(K)-C^\star)/\gamma^\star$, so $C(K)-C(K')\ge(2-\chi)(\varkappa/\gamma^\star)(C(K)-C^\star)$, and subtracting this bound from $C(K)-C^\star$ gives~\eqref{eq:gnrate}.
For the bounds on $(2-\chi)\varkappa$, multiplying $MSM\preceq\chi M$ on both sides by $M^{-1/2}$ gives $M^{1/2}SM^{1/2}\preceq\chi I$.
The matrices $M^{1/2}SM^{1/2}$ and $S^{1/2}MS^{1/2}$ have the same eigenvalues, and $S^{1/2}MS^{1/2}\succ0$, so $0<\varkappa\le\chi$.
Hence $0<(2-\chi)\varkappa\le(2-\chi)\chi\le1\le\gamma^\star$, where $(2-\chi)\chi\le1$ because $1-(2-\chi)\chi=(1-\chi)^2$.
\end{proof}

The proofs of Theorems~\ref{thm:mb} and~\ref{thm:adaptive} need bounds that hold for every gain along the iterations and a compact set that contains these gains, and the next lemma provides both.

\begin{applemma}[Bounds]\label{lem:bounds}
For $K\in\cS$, $\|P_K\|\le C(K)/\usig(\Sigma_W)$, $\|\Sigma_K\|\le C(K)/\usig(Q)$ and $\|K\|^2\le C(K)/(\usig(R)\usig(\Sigma_W))$, and $\cS(v)$ is compact.
\end{applemma}

\begin{proof}
The first two bounds are~\cite[Lemma~13]{fazel2018global}, and the third extends~\cite[Lemma~10]{zhao2025pgac} to a general $\Sigma_W$.
The bounds follow from $C(K)=\Tr(P_K\Sigma_W)=\Tr((Q+K^\top RK)\Sigma_K)$ with $P_K\succeq0$ and $\Sigma_K\succeq\Sigma_W\succeq\usig(\Sigma_W)I_n$; in particular $\cS(v)$ is bounded.
For closedness, let $K_j\in\cS(v)$ converge to $\bar K$; since $\Sigma_{K_j}=\sum_{k\ge0}A_{K_j}^k\Sigma_W(A_{K_j}^\top)^k$, the partial sums $\sum_{k\le N}A_{K_j}^k(A_{K_j}^\top)^k\preceq\Sigma_{K_j}/\usig(\Sigma_W)\preceq\frac{v}{\usig(Q)\usig(\Sigma_W)}I$ pass to the limit for every $N$, so $A_{\bar K}$ is Schur, and $C(\bar K)\le v$ by continuity of $C$ on $\cS$.
\end{proof}

The last lemma of this appendix bounds the change of the solution of a Lyapunov equation when its Schur matrix $L$ and its constant term $Y$ change, and Appendix~\ref{app:adaptive} uses it to prove Facts~\ref{f:pert} and~\ref{f:lip}.
It is a variant of~\cite[Lemma~11]{zhao2025pgac}, which treats $\Sigma=I_n+L\Sigma L^\top$ and gives a bound of this kind with other constants.
The general weight $Y$ is needed for the equation of $P_K$, whose weight is $Q+K^\top RK$, and for the equation of $\Sigma_K$ with a general $\Sigma_W$; the proof of Fact~\ref{f:lip} uses the change of the weight.

\begin{applemma}[Lyapunov perturbation]\label{lem:lyap}
Let $L$ be Schur, $Y\succeq yI$ with $y>0$, and let $P$ solve $P=Y+L^\top PL$.
Let $L'$ and a symmetric $Y'$ be such that $\delta:=\|L'-L\|$ satisfies $\beta:=\|P\|\delta(2\|L\|+\delta)\le\frac y2$.
Then $L'$ is Schur, and the solution $P'$ of $P'=Y'+L'^\top P'L'$ satisfies
\begin{equation}\label{eq:lyap}
\|P'-P\|\le\tfrac{2\|P\|}{y}\big(\|Y'-Y\|+\beta\big).
\end{equation}
The same holds with $L$ and $L'$ replaced by $L^\top$ and $L'^\top$, which covers the equations $\Sigma=Y+L\Sigma L^\top$.
\end{applemma}

\begin{proof}
The idea is that $P$ certifies the stability of $L'$ with margin $\frac y2$, and summing the Lyapunov difference along $L'$ then bounds $P'-P$.
The matrix $\Delta:=L'^\top PL'-L^\top PL=(L'-L)^\top PL'+L^\top P(L'-L)$ satisfies $\|\Delta\|\le\|P\|\delta(\|L'\|+\|L\|)\le\beta$, so $P-L'^\top PL'=Y-\Delta\succeq\frac y2I$.
For an eigenvector $v$ of $L'$ with eigenvalue $\zeta_v$, the bound $P-L'^\top PL'\succeq\frac y2I$ gives $(1-|\zeta_v|^2)v^\ast Pv>0$ with $P\succeq Y\succ0$, so $L'$ is Schur.
Summing $L'^{k\top}(P-L'^\top PL')L'^k$ over $0\le k\le N$ telescopes to $P-L'^{(N+1)\top}PL'^{N+1}\preceq P$, so $G:=\sum_{k\ge0}L'^{k\top}L'^k\preceq\frac2yP$.
Subtracting the two Lyapunov equations gives $P'-P=Y'-Y+\Delta+L'^\top(P'-P)L'$, so $P'-P=\sum_{k\ge0}L'^{k\top}(Y'-Y+\Delta)L'^k$, whose norm is at most $\|G\|\,\|Y'-Y+\Delta\|\le\frac{2\|P\|}y(\|Y'-Y\|+\beta)$, which is~\eqref{eq:lyap}.
The transposed case is the same argument applied to $L^\top$ and $L'^\top$, which have the norms of $L$ and $L'$.
\end{proof}

\section{Proofs of Sections~\ref{sec:kl} and~\ref{sec:alg}}\label{app:kl}

The proofs rest on the Gibbs variational identity, applied state by state.

\begin{applemma}[{Gibbs variational identity~\cite[Thm.~1 and Lemma~1]{theodorou2015nonlinear}}]\label{lem:gibbs}
Let $\tau>0$, $\mu$ a distribution and $f$ a function with $\mathcal Z:=\EE_{Y\sim\mu}[e^{-f(Y)/\tau}]<\infty$, and $\mu^\star:=\mu e^{-f/\tau}/\mathcal Z$.
For every $\mu'$ with $\KL(\mu'\|\mu)<\infty$ and $\EE_{Y\sim\mu'}|f(Y)|<\infty$,
\begin{equation}\label{eq:gibbs}
\EE_{Y\sim\mu'}[f(Y)]+\tau\KL(\mu'\|\mu)=-\tau\log\mathcal Z+\tau\KL(\mu'\|\mu^\star),
\end{equation}
so $\mu^\star$ is the unique minimizer of the left side.
\end{applemma}

The cited results state the variational form of~\eqref{eq:gibbs} and its minimizer; the identity itself follows by inserting the density of $\mu^\star$ into $\KL(\mu'\|\mu^\star)$.

\runin{Proof of Proposition~\ref{prop:soft}}
The proof evaluates $\pi$ through its Bellman equation for~(a), computes the Gibbs minimizer by completing a square for~(b), shows that the new gain is stabilizing and compares the two policies by telescoping for~(c), and then specializes to $K_\rho=K$ for~(d).
Throughout, $P:=P^\rho_K$, $S:=S^\rho_K$, $N:=B^\top PA$ and $E:=E^\rho_K=SK+N+D(K-K_\rho)$.
We first prove~(a).
Since $\pi(\cdot\mid x)$ and $\rho(\cdot\mid x)$ are Gaussian with means $Kx$ and $K_\rho x$, the penalty is $\tau\kappa^\pi_\rho(x)=\varpi+x^\top(K-K_\rho)^\top D(K-K_\rho)x$ with $\varpi:=\tau\KL(\Nor(0,\Psi)\|\Nor(0,\Psi_\rho))$.
With $V(x):=x^\top Px$, $\Qf(x,u):=c(x,u)+\EE_{X'\sim\Nor(Ax+Bu,\Sigma_W)}[V(X')]=x^\top Qx+u^\top Ru+(Ax+Bu)^\top P(Ax+Bu)+\Tr(P\Sigma_W)$.
Averaging over $U=Kx+\Omega$, $\Omega\sim\Nor(0,\Psi)$, removes the term linear in $\Omega$ and gives
\begin{equation}\label{eq:bellpi}
\begin{aligned}
&\EE_{U\sim\pi(\cdot\mid x)}[\Qf(x,U)]+\tau\kappa^\pi_\rho(x)\\
&\quad=V(x)+\Tr(P\Sigma_W)+\Tr(S\Psi)+\varpi,
\end{aligned}
\end{equation}
because $Q+K^\top RK+(K-K_\rho)^\top D(K-K_\rho)+A_K^\top PA_K=P$.
Hence $V$ and $\bar J:=\Tr(P\Sigma_W)+\Tr(S\Psi)+\varpi$ solve the Bellman equation of $\pi$ for (P).
Summing~\eqref{eq:bellpi} at $X_t$ along $\PP^\pi$ telescopes the terms in $V$, and since $A_K$ is Schur and the inputs are Gaussian, $\EE_{\PP^\pi}[V(X_T)]$ is bounded in $T$, so $J_\rho(\pi)=\bar J$ and $V_\pi=V$ up to a constant.
Next we prove~(b) by computing the minimizer of~\eqref{eq:softpi}.
We apply Lemma~\ref{lem:gibbs} with $\mu=\rho(\cdot\mid x)$ and $f=\Qf(x,\cdot)$; the function $f$ is bounded below since $S\succ0$, so $\mathcal Z<\infty$.
Since $V_\pi=V$ up to a constant, $\Qf_\pi$ and $\Qf$ differ by a constant, so the left side of~\eqref{eq:gibbs} with these choices equals the objective of~\eqref{eq:softpi} up to a constant, and the minimizer of~\eqref{eq:softpi} is the distribution $\mu^\star$ of Lemma~\ref{lem:gibbs}, proportional to $\rho(u\mid x)e^{-\Qf(x,u)/\tau}$.
Its exponent in $u$ is
\begin{equation}\label{eq:gibbsexp}
-\tfrac12u^\top\big(\Psi_\rho^{-1}+\tfrac2\tau S\big)u+u^\top\big(\Psi_\rho^{-1}K_\rho-\tfrac2\tau N\big)x+\text{terms in }x,
\end{equation}
so $\pi'(\cdot\mid x)$ is Gaussian with precision $\Psi_\rho^{-1}+\frac2\tau S=\frac2\tau(D+S)=\frac2\tau H$ and mean $H^{-1}(DK_\rho-N)x=H^{-1}(HK-E)x$, which is~\eqref{eq:softgen}.
We now prove~(c).
Write $c_\rho(x,u):=c(x,u)+(u-K_\rho x)^\top D(u-K_\rho x)$.
For this stage cost, $P$ solves the Lyapunov equation of the gain $K$, and expanding as in~\eqref{eq:defect} gives, for every $K'$ with $\Delta:=K'-K$,
\begin{equation*}
\begin{aligned}
&Q+K'^\top RK'+(K'-K_\rho)^\top D(K'-K_\rho)+A_{K'}^\top PA_{K'}\\
&\quad-P=\Delta^\top E+E^\top\Delta+\Delta^\top H\Delta.
\end{aligned}
\end{equation*}
By~\eqref{eq:softgen}, $\Delta=-H^{-1}E$, so the right side of this expansion is $-2E^\top H^{-1}E+E^\top H^{-1}E=-E^\top H^{-1}E\preceq0$.
The matrix $P$ solves the Lyapunov equation of the Schur matrix $A_K$ with a constant term at least $Q$, so $P\succeq Q$.
Hence $A_{K'}^\top PA_{K'}-P\preceq-Q\prec0$ with $P\succeq Q\succ0$, so $A_{K'}$ is Schur by Lyapunov's stability theorem, that is, $K'\in\cS$.
By the optimality of $\pi'$ in~\eqref{eq:softpi} and by~\eqref{eq:bellpi}, for every $x$,
\begin{equation}\label{eq:improve}
\EE_{U\sim\pi'(\cdot\mid x)}[\Qf(x,U)]+\tau\kappa^{\pi'}_\rho(x)\le V(x)+\bar J.
\end{equation}
Evaluating~\eqref{eq:improve} at $X_t$ under $\PP^{\pi'}$, using $\EE_{\PP^{\pi'}}[\Qf(X_t,U_t)]=\EE_{\PP^{\pi'}}[c(X_t,U_t)+V(X_{t+1})]$ and summing over $t<T$ gives $\EE_{\PP^{\pi'}}[\sum_{t<T}(c(X_t,U_t)+\tau\kappa^{\pi'}_\rho(X_t))]\le T\bar J+\EE_{\PP^{\pi'}}[V(X_0)-V(X_T)]\le T\bar J+\EE_{\PP^{\pi'}}[V(X_0)]$, since $V\ge0$; the expectations are finite because $K'\in\cS$ and $\kappa^{\pi'}_\rho$ is quadratic in $x$.
Dividing by $T$ gives $J_\rho(\pi')\le\bar J=J_\rho(\pi)$.
Finally, we prove~(d).
If $K_\rho=K$, the Lyapunov equation of $P^\rho_K$ is that of $P_K$, so $S=S_K$, $E=E_K$, $\Qf=\Qf_K+C(K)$ by~\eqref{eq:Q}, and~\eqref{eq:softgen} becomes~\eqref{eq:softclosed}.
Since $H\succeq S_K$, $M:=H^{-1}$ satisfies $MS_KM\preceq M$, and the gain of~\eqref{eq:softclosed} is $K'=K-ME_K$.
Lemma~\ref{lem:loewner} with this $M$ and $\chi=1$ therefore gives $P_{K'}\preceq P_K-E_K^\top ME_K=P_K-E_K^\top H^{-1}E_K$.\hfill$\square$

\runin{Proof of Proposition~\ref{prop:price}}
The closed loop $X_{t+1}=A_KX_t+B\Omega_t+W_t$ is stable with the sum $B\Omega_t+W_t$ of covariance $\Sigma_W+B\Psi B^\top$, so by the duality of the two Lyapunov equations its average cost is $\Tr(P_K(\Sigma_W+B\Psi B^\top))+\Tr(R\Psi)=C(K)+\Tr(S_K\Psi)$.
With $\lambda_1,\dots,\lambda_m$ the eigenvalues of $S_K^{1/2}\Psi S_K^{1/2}$, $\Tr(S_K\Psi)=\sum_i\lambda_i\ge m(\prod_i\lambda_i)^{1/m}$ by the inequality of arithmetic and geometric means, with equality exactly when the $\lambda_i$ are equal.\hfill$\square$

\runin{Proof of Proposition~\ref{prop:corr}}
For (i), $H=S_K+D$ by~\eqref{eq:softclosed}, and $D=\frac\tau2\Psi_\rho^{-1}$ is a bijection between $\Psi_\rho\succ0$ and $D\succ0$, so $H\succ S_K$ is reached exactly by $D=H-S_K$.
For (ii), $H=\tfrac{S_K}{2\eta}$ gives $D=\frac{1-2\eta}{2\eta}S_K\succ0$ for $\eta<\frac12$, hence $\Psi_\rho=\frac\tau2D^{-1}=\frac{\eta\tau}{1-2\eta}S_K^{-1}$ and $\Psi'=\frac\tau2H^{-1}=\eta\tau S_K^{-1}$.
For (iii), $H=\tfrac{1}{2\eta}I$ gives $D=\frac1{2\eta}I-S_K\succ0$ for $\eta<\tfrac{1}{2\lambda_{\max}(S_K)}$, so $\Psi_\rho=\tau(\eta^{-1}I-2S_K)^{-1}$ and $\Psi'=\eta\tau I$.
For (iv), the first claim is continuity of~\eqref{eq:softclosed} in $D$.
If a centred reference with scaled precision $D\succ0$ gave the step $K-2\eta S_K^{-1}E_K$, then for $x$ with $E_Kx\neq0$ the vector $v:=2\eta S_K^{-1}E_Kx\neq0$ would satisfy $(S_K+D)v=E_Kx=\frac1{2\eta}S_Kv$, so $v^\top Dv=\frac{1-2\eta}{2\eta}v^\top S_Kv\le0$ for $\eta\ge\frac12$, which contradicts $v^\top Dv>0$.\hfill$\square$

Lemma~\ref{lem:mm} gives the fixed objective of Section~\ref{sec:alg}: no iteration of KLPI increases $J_\vartheta$ of~\eqref{eq:Jent}.
The majorization--minimization (MM) principle~\cite{lange2016mm} decreases an objective $F$ through a surrogate that lies above $F$ and touches it at the current iterate, and a step that decreases the surrogate without minimizing it also does not increase $F$.

\begin{applemma}[Majorization and fixed objective]\label{lem:mm}
Let $\tau_t=\tau$, $\alpha\in(0,1]$, $\vartheta:=(1-\alpha)\tau$, and iterate~\eqref{eq:alg} with the true model from $K_0\in\cS$, $H_0\succ0$, and let $\pi_t:=\Nor(K_tx,\frac\tau2H_t^{-1})$.
For linear-Gaussian $\pi$ as in~\eqref{eq:Jent}, with $d_\pi$ the stationary state distribution of the closed loop under $\pi$, let
\begin{equation}\label{eq:surrogate}
G_t(\pi):=J_\vartheta(\pi)+\alpha\tau\,\EE_{X\sim d_\pi}\big[\kappa^\pi_{\pi_t}(X)\big].
\end{equation}
\pitem{(a)} $G_t=J_{\rho_t}-\tau\log Z_t$ with a constant $Z_t$ that does not depend on $\pi$, $G_t\ge J_\vartheta$ with equality at $\pi_t$, and
\begin{equation}\label{eq:mm}
J_\vartheta(\pi_{t+1})\le G_t(\pi_{t+1})\le G_t(\pi_t)=J_\vartheta(\pi_t).
\end{equation}
\pitem{(b)} In every state $x$, $\pi_{t+1}(\cdot\mid x)$ minimizes $\EE_{U\sim\mu}[\Qf_{\pi_t}(x,U)]-\vartheta h(\mu)+\alpha\tau\KL(\mu\|\pi_t(\cdot\mid x))$ over distributions $\mu$, which is the form of a policy mirror-descent step on the entropy-regularized cost~\cite{lan2023policy,cen2022fast}.
\pitem{(c)} For $\vartheta>0$, among the linear-Gaussian policies of~\eqref{eq:Jent}, $J_\vartheta$ has the unique minimizer $\Nor\big(K^\star x,\frac\vartheta2(S^\star)^{-1}\big)$.
\end{applemma}

Each iteration of KLPI is thus a generalized MM step: one policy improvement does not increase the surrogate $G_t$, although it need not minimize it.
For $\alpha=1$, $J_0=C(K,\Psi)$, and~\eqref{eq:mm} is the average-cost counterpart of the finite-horizon majorization of~\cite{lefebvre2024probabilistic}.
The infimum $C^\star$ of $J_0$ is approached only as $\Psi\to0$, so (c) requires $\vartheta>0$.

\runin{Proof of Lemma~\ref{lem:mm}}
We first express $G_t$ through the objective $J_{\rho_t}$ of the policy improvement, then deduce~\eqref{eq:mm} and~(b) from Proposition~\ref{prop:soft}, and finally minimize $J_\vartheta$ in closed form.
For a density $\mu$ on $\mathbb R^m$ and $x$ fixed, write $\pi_t=\pi_t(\cdot\mid x)$ and $Z_t:=\int\pi_t(u)^\alpha du$, so that $\rho_t=\pi_t^\alpha/Z_t$; for $\pi_t=\Nor(K_tx,\Psi_t)$, $Z_t$ depends on $\Psi_t$ and $\alpha$ only.
Since $\log(\mu/\rho_t)=\alpha\log(\mu/\pi_t)+(1-\alpha)\log\mu+\log Z_t$, averaging over $U\sim\mu$ gives
\begin{equation}\label{eq:klsplit}
\tau\KL(\mu\|\rho_t)=\alpha\tau\KL(\mu\|\pi_t)-\vartheta h(\mu)+\tau\log Z_t.
\end{equation}
For $\pi=\Nor(Kx,\Psi)$ with $K\in\cS$, $h(\pi(\cdot\mid x))=h(\Nor(0,\Psi))$, which is $\frac12\log\det\Psi$ plus a constant that depends on $m$ only, for every $x$, and $\kappa^\pi_{\pi_t}(x)$ is a quadratic function of $x$ plus a constant.
Under $\PP^\pi$ the second moment of $X_t$ converges to that of $d_\pi=\Nor(0,\Sigma_{K,\Psi})$, $\Sigma_{K,\Psi}=\Sigma_W+B\Psi B^\top+A_K\Sigma_{K,\Psi}A_K^\top$, so every time average in~\eqref{eq:P} converges, and by~\eqref{eq:klsplit} and Proposition~\ref{prop:price}, $J_{\rho_t}(\pi)=G_t(\pi)+\tau\log Z_t$.
Since $\kappa^\pi_{\pi_t}\ge0$, $G_t\ge J_\vartheta$, with equality at $\pi=\pi_t$.
Next we prove~\eqref{eq:mm}.
By Proposition~\ref{prop:alg} with $\alpha_t=\alpha$, $\pi_{t+1}$ is the policy improvement of $\pi_t$ against $\rho_t$, and Proposition~\ref{prop:soft}(c), applied with $\pi=\pi_t$ and $\rho=\rho_t$, both centred at $K_t$, gives $K_{t+1}\in\cS$ and $J_{\rho_t}(\pi_{t+1})\le J_{\rho_t}(\pi_t)$.
Subtracting $\tau\log Z_t$ from both sides gives $G_t(\pi_{t+1})\le G_t(\pi_t)$, and the bound $G_t\ge J_\vartheta$ with equality at $\pi_t$ gives the two remaining relations of~\eqref{eq:mm}.
For (b), by~\eqref{eq:klsplit} the objective of (b) equals that of~\eqref{eq:softpi} with the reference $\rho_t$ minus the constant $\tau\log Z_t$, and the state-action value $\Qf_{\pi_t}$ of~\eqref{eq:softpi} is $\Qf_{K_t}+C(K_t)$ by Proposition~\ref{prop:soft}(d), whatever the reference centred at $K_t$.
Finally, we prove~(c).
By Proposition~\ref{prop:price}, $J_\vartheta=C(K)+\Tr(S_K\Psi)-\frac\vartheta2\log\det\Psi+\text{const}$, which is strictly convex in $\Psi\succ0$ with the unique minimizer $\Psi=\frac\vartheta2S_K^{-1}$ and minimum $C(K)+\frac\vartheta2\log\det S_K+\text{const}$.
Both terms of this minimum are minimized at $K^\star$, the first uniquely by Lemma~\ref{lem:costdiff}, and the second because $S_K=R+B^\top P_KB\succeq S^\star$, since $P_K\succeq P^\star$ for $K\in\cS$ by Lemma~\ref{lem:costdiff}.\hfill$\square$

\section{Proof of Theorem~\ref{thm:mb}}\label{app:mb}

The step is an instance of Lemma~\ref{lem:gn} with the preconditioner $M=H_{g+1}^{-1}$, and the work lies in bounding the accumulated Hessian $H_{g+1}$, which sets the fraction of the gap removed per step.
We first show that every step decreases the cost, then bound the accumulated Hessian, and finally treat the cases $\alpha<1$ and $\alpha=1$.
If $K_g\in\cS$ and $H_g\succeq0$, then $H_{g+1}=\alpha H_g+S_{K_g}\succeq S_{K_g}$, so $M:=H_{g+1}^{-1}$ satisfies $MS_{K_g}M\preceq M$, and Lemmas~\ref{lem:loewner} and~\ref{lem:gn} with $\chi=1$ give $K_{g+1}\in\cS$, the Loewner decrease, and
\begin{equation}\label{eq:mbstep}
C(K_{g+1})-C^\star\le\Big(1-\frac{\underline{s}}{\lambda_{\max}(H_{g+1})\gamma^\star}\Big)\big(C(K_g)-C^\star\big),
\end{equation}
since $\varkappa\ge\lambda_{\min}(M)\underline{s}=\underline{s}/\lambda_{\max}(H_{g+1})$.
The regularization weights do not enter.
Since $K_0\in\cS$, $H_0\succ0$ and $H_{g+1}\succeq S_{K_g}\succ0$, induction on $g$ shows that $K_g\in\cS$ and $P_{K_{g+1}}\preceq P_{K_g}$ for every $g$.
Next we bound the accumulated Hessian.
By $P_{K_k}\preceq P_{K_0}$ and Lemma~\ref{lem:bounds}, $\|P_{K_k}\|\le\|P_{K_0}\|\le C(K_0)/\usig(\Sigma_W)$, so $R\preceq S_{K_k}\preceq\bar sI$, and unrolling the recursion,
\begin{equation}\label{eq:Hbound}
\lambda_{\max}(H_{g+1})\le\alpha^{g+1}\|H_0\|+\bar s\textstyle\sum_{k\le g}\alpha^{g-k}.
\end{equation}
We now treat the two cases.
For $\alpha<1$,~\eqref{eq:Hbound} gives $H_g\preceq(\|H_0\|+\tfrac{\bar s}{1-\alpha})I$ for all $g$, hence the floor on $\Psi_g=\frac{\tau_g}2H_g^{-1}$, and~\eqref{eq:mbstep} gives the factor $1-\theta_\alpha$, since $\theta_\alpha=\underline{s}/((\|H_0\|+\frac{\bar s}{1-\alpha})\gamma^\star)\le \underline{s}/(\lambda_{\max}(H_{g+1})\gamma^\star)$.
Moreover, $\underline{s}/(\lambda_{\max}(H_{g+1})\gamma^\star)\le\varkappa/\gamma^\star\le1$, by the bound on $\varkappa$ stated after~\eqref{eq:mbstep} and by Lemma~\ref{lem:gn} with $\chi=1$; hence $\theta_\alpha\le1$.
The gains stay in the compact set $\cS(C(K_0))$ (Lemma~\ref{lem:bounds}) and their costs tend to $C^\star$; since $K^\star$ is the only gain of cost $C^\star$ (Lemma~\ref{lem:costdiff}) and $C$ is continuous on $\cS$, $K_g\to K^\star$, so $S_{K_g}\to S^\star$, and the affine recursion gives $H_g\to \tfrac{S^\star}{1-\alpha}$.
For $\alpha=1$,~\eqref{eq:Hbound} gives $\lambda_{\max}(H_{g+1})\le\|H_0\|+(g+1)\bar s=\bar s(c_1+g)$, so the factor in~\eqref{eq:mbstep} is at most $1-c_2/(c_1+g)$, with $c_2\le1\le c_1$ because $\gamma^\star\ge1$ and $\underline{s}\le\bar s$, the latter since $Q\preceq P_{K_0}$.
Then
\begin{equation}\label{eq:sublin}
\prod_{k<g}\Big(1-\frac{c_2}{c_1+k}\Big)\le\exp\Big(-c_2\sum_{k<g}\frac1{c_1+k}\Big)\le\Big(\frac{c_1}{c_1+g}\Big)^{c_2},
\end{equation}
since $\sum_{k<g}\frac1{c_1+k}\ge\int_0^g\frac{dy}{c_1+y}=\log\frac{c_1+g}{c_1}$.
Finally $R\preceq S_{K_k}\preceq\bar sI$ gives $(\lambda_{\min}(H_0)+g\lambda_{\min}(R))I\preceq H_g\preceq(\|H_0\|+g\bar s)I$, which bounds $\Psi_g=\frac{\tau_g}2H_g^{-1}$ on both sides.\hfill$\square$

\section{Proofs of Section~\ref{sec:adaptive} and Proposition~\ref{prop:regret}}\label{app:adaptive}

The proof of Theorem~\ref{thm:adaptive} follows the architecture of the proofs of~\cite[Theorems~1 and~3]{zhao2025pgac} and uses four facts, whose counterparts for $\Sigma_W=I_n$ are in~\cite{zhao2025pgac}.

\runin{Stability notions~\cite{cohen2019learning,zhao2025pgac}}
A gain $K$ is $(b,\varepsilon)$-strongly stable, $b\ge1$, $\varepsilon\in(0,1]$, if $A_K=\mathcal GL\mathcal G^{-1}$ with $\mathcal G\succ0$, $\|L\|\le1-\varepsilon$, $\|\mathcal G\|\le b$, $\|\mathcal G^{-1}\|\le1$ and $\|K\|\le b$.
A sequence $(K_t)$ is sequentially stable if every $K_t$ is $(b,\varepsilon)$-strongly stable with common constants through matrices $\mathcal G_t$ with $\|\mathcal G_{t+1}^{-1}\mathcal G_t\|\le1+\varepsilon/2$.

\runin{Four facts}
Let $\bar c:=C(K_{t_0})+2$.
For $K\in\cS(\bar c)$, the numbers $\ell_P:=\bar c/\usig(\Sigma_W)$, $\ell_\Sigma:=\bar c/\usig(Q)$, $\ell_K:=(\bar c/(\usig(R)\usig(\Sigma_W)))^{1/2}$ and $\ell_A:=(\ell_P/\usig(Q))^{1/2}$ bound $\|P_K\|$, $\|\Sigma_K\|$, $\|K\|$ and $\|A_K\|$, the first three by Lemma~\ref{lem:bounds} and $\|A_K\|\le\ell_A$ because $\usig(Q)\|A_Kv\|^2\le v^\top A_K^\top P_KA_Kv\le v^\top P_Kv$ for every $v$.
They also satisfy $\ell_P\ge\|P_K\|\ge\usig(Q)$ and $\ell_\Sigma\ge\|\Sigma_K\|\ge\usig(\Sigma_W)$.
There are constants $p_1,p_2,p_4,p_5,p_6,p_9>0$, numbered as their counterparts in~\cite{zhao2025pgac} and given in the proofs that follow the facts, depending only on $(A,B,Q,R,\Sigma_W,\bar c)$, such that for every $K\in\cS(\bar c)$ the following hold, where, for each $t\ge t_0$, a hat denotes an object computed on the plant estimate $(\hat A_t,\hat B_t)$ (Fact~\ref{f:pert} is applied at index $t+1$ to $\hat S_t$ and $\hat E_t$ of Algorithm~\ref{alg:akpi}, and at index $t_0$ to $\hat S_{\mathrm{off}}$).

\begin{fact}[Sensitivity to the plant estimate]\label{f:pert}
Under Assumption~\ref{ass:pe}, if $\SNR_t\ge1/p_2$, then $\hat A_t+\hat B_tK$ is Schur, $\|\hat P_K-P_K\|\le p_1/\SNR_t$ and $\|\hat E_K-E_K\|\le p_9/\SNR_t$.
\end{fact}
\begin{fact}[Local Lipschitz continuity]\label{f:lip}
If $\|\tilde K-K\|\le p_4$, then $\tilde K\in\cS$, $\|\Sigma_{\tilde K}-\Sigma_K\|\le p_5\|\tilde K-K\|$ and $|C(\tilde K)-C(K)|\le p_6\|\tilde K-K\|$.
\end{fact}
\begin{fact}[Uniform strong stability]\label{f:strong}
$K$ is $(\nu_3,\nu_4)$-strongly stable with $\mathcal G=(\Sigma_K/\usig(\Sigma_W))^{1/2}$, $\nu_3:=\big(\bar c/(\usig(\Sigma_W)\min(\usig(Q),\usig(R)))\big)^{1/2}$ and $\nu_4:=1-(1-\nu_3^{-2})^{1/2}$.
\end{fact}
\begin{fact}[State bound]\label{f:state}
If $K_t\in\cS(\bar c)$ and $\|\Sigma_{K_{t+1}}^{-1}\Sigma_{K_t}\|\le1+\nu_4/2$ for all $t\ge t_0$, then $(K_t)$ is sequentially stable and~\eqref{eq:statebound} holds for $X_{t+1}=A_{K_t}X_t+B\Omega_t+W_t$.
\end{fact}

\runin{Proof of Fact~\ref{f:pert}}
Let $\xi:=1/\SNR_t$, $L:=A_K$ and $L':=\hat A_t+\hat B_tK$.
By~\cite[Lemma~1]{zhao2025pgac}, $\delta:=\|L'-L\|=\|[\hat B_t-B\ \ \hat A_t-A]\big[\begin{smallmatrix}K\\ I_n\end{smallmatrix}\big]\|\le(1+\ell_K)\xi$.
Let $p_2:=\usig(Q)/(2\ell_P(1+\ell_K)(2\ell_A+1))$ and $p_1:=2\ell_P^2(1+\ell_K)(2\ell_A+1)/\usig(Q)$, so that $p_1p_2=\ell_P$.
If $\xi\le p_2$, then $\delta\le\usig(Q)/(2\ell_P(2\ell_A+1))\le\frac12$, so Lemma~\ref{lem:lyap} with $P=P_K$ and $Y=Y'=Q+K^\top RK\succeq\usig(Q)I$ applies with $\beta\le\ell_P\delta(2\ell_A+1)\le\frac{\usig(Q)}2$.
Lemma~\ref{lem:lyap} shows that $L'$ is Schur and gives $\|\hat P_K-P_K\|\le\frac{2\ell_P}{\usig(Q)}\ell_P(2\ell_A+1)(1+\ell_K)\xi=p_1\xi$.
Since $E_K=RK+B^\top P_KA_K$ and $\hat E_K=RK+\hat B_t^\top\hat P_KL'$,
\begin{equation*}
\begin{aligned}
\hat E_K-E_K={}&(\hat B_t-B)^\top\hat P_KL'+B^\top(\hat P_K-P_K)L'\\
&+B^\top P_K(L'-L),
\end{aligned}
\end{equation*}
and with $\|\hat B_t-B\|\le\xi$, $\|\hat P_K\|\le\ell_P+p_1p_2=2\ell_P$ and $\|L'\|\le\ell_A+1$, this gives $\|\hat E_K-E_K\|\le p_9\xi$ with $p_9:=(\ell_A+1)(2\ell_P+\|B\|p_1)+\|B\|\ell_P(1+\ell_K)$.
For $\Sigma_W=I_n$, counterparts of the bounds on $\hat P_K$ and $\hat E_K$ appear in~\cite[Lemma~8]{zhao2025deepo} and~\cite[Lemmas~12 and~23]{zhao2025pgac}.

\runin{Proof of Fact~\ref{f:lip}}
Let $d:=\|\tilde K-K\|\le p_4$ with $p_4:=\min\{1,\,\usig(Q)/(2\ell_P\|B\|(2\ell_A+1)),\,\usig(\Sigma_W)/(2\ell_\Sigma\|B\|(2\ell_A+1))\}$.
Then $\delta:=\|A_{\tilde K}-A_K\|\le\|B\|d\le\frac12$, and $\|\tilde K^\top R\tilde K-K^\top RK\|\le\|R\|d(\|\tilde K\|+\|K\|)\le\|R\|(2\ell_K+1)d$.
Lemma~\ref{lem:lyap} with $L=A_K$, $P=P_K$, $Y=Q+K^\top RK$ and $Y'=Q+\tilde K^\top R\tilde K$ applies with $\beta\le\ell_P\|B\|(2\ell_A+1)d\le\frac{\usig(Q)}2$, so $\tilde K\in\cS$ and
\begin{equation*}
\|P_{\tilde K}-P_K\|\le\tfrac{2\ell_P}{\usig(Q)}\big(\|R\|(2\ell_K+1)+\ell_P\|B\|(2\ell_A+1)\big)d.
\end{equation*}
Since $|C(\tilde K)-C(K)|=|\Tr((P_{\tilde K}-P_K)\Sigma_W)|\le\Tr(\Sigma_W)\|P_{\tilde K}-P_K\|$, this gives the cost bound with $p_6:=\frac{2\ell_P\Tr(\Sigma_W)}{\usig(Q)}\big(\|R\|(2\ell_K+1)+\ell_P\|B\|(2\ell_A+1)\big)$.
The transposed form of Lemma~\ref{lem:lyap} with $P=\Sigma_K$ and $Y=Y'=\Sigma_W$ applies with $\beta\le\ell_\Sigma\|B\|(2\ell_A+1)d\le\frac{\usig(\Sigma_W)}2$ and gives the covariance bound with $p_5:=2\ell_\Sigma^2\|B\|(2\ell_A+1)/\usig(\Sigma_W)$.
For $\Sigma_W=I_n$, the counterpart of Fact~\ref{f:lip} is~\cite[Lemma~13]{zhao2025pgac}.

\runin{Proofs of Facts~\ref{f:strong} and~\ref{f:state}}
Fact~\ref{f:strong} repeats the argument of~\cite{cohen2019learning,zhao2025pgac}: by the Lyapunov equation, $L:=\Sigma_K^{-1/2}A_K\Sigma_K^{1/2}$ satisfies $LL^\top=I_n-\Sigma_K^{-1/2}\Sigma_W\Sigma_K^{-1/2}\preceq(1-\usig(\Sigma_W)/\|\Sigma_K\|)I_n$, so Lemma~\ref{lem:bounds} gives $\|L\|^2\le1-\nu_3^{-2}$, that is, $\|L\|\le1-\nu_4$, while $\Sigma_K\succeq\usig(\Sigma_W)I_n$ gives $\|\mathcal G^{-1}\|\le1$, and Lemma~\ref{lem:bounds} gives $\|\mathcal G\|\le \nu_3$ and $\|K\|\le \nu_3$, and $\nu_3\ge1$ because $\bar c\ge C(K)\ge\Tr(Q\Sigma_W)\ge\usig(Q)\usig(\Sigma_W)$.
Fact~\ref{f:state} uses $\mathcal G_t=(\Sigma_{K_t}/\usig(\Sigma_W))^{1/2}$ and $\|\mathcal G_{t+1}^{-1}\mathcal G_t\|^2=\lambda_{\max}(\Sigma_{K_{t+1}}^{-1/2}\Sigma_{K_t}\Sigma_{K_{t+1}}^{-1/2})\le\|\Sigma_{K_{t+1}}^{-1}\Sigma_{K_t}\|$, with $\sqrt{1+x}\le1+x$; the state bound is then the argument of~\cite[Lemma~3]{cohen2019learning} reproduced in~\cite{zhao2025pgac}.

The proof of Theorem~\ref{thm:adaptive} is an induction on $t$ that keeps the gain in $\cS(\bar c)$ and $H_t$ in the band of (C3).
Lemma~\ref{lem:hessband} bounds the estimated input Hessian of every gain in $\cS(\bar c)$, and Lemma~\ref{lem:onestep} uses this bound to carry out one step of the induction.

\runin{Constants}
Let $\bar p:=\bar c/\usig(\Sigma_W)+p_1p_2$, $\bar b:=\|B\|+p_2$, $\sa:=\|R\|+\bar b^2\bar p$, $r_0:=\underline{s}/(2\sa\gamma^\star)$, $p_S:=\bar b(\bar p+\|B\|p_1)+\|B\|\bar c/\usig(\Sigma_W)$, $\nu_5:=2p_6p_9/\underline{s}$, $\nu_1:=\max\{1/p_2,\,8p_S/\underline{s},\,8p_9/(\underline{s}p_4),\,2\nu_5/r_0\}$, $p_8:=\sup_{K\in\cS(\bar c)}\|E_K\|+p_9p_2$, finite by Lemma~\ref{lem:bounds}, and $\nu_2:=\frac{\underline{s}}2\min\{p_4/p_8,\ \usig(\Sigma_W)\nu_4/(4\nu_3^2p_5p_8)\}$.
Since $Q\preceq P_K$ and $\|P_K\|\le\bar c/\usig(\Sigma_W)$ on $\cS(\bar c)$ (Lemma~\ref{lem:bounds}), $\underline{s}\le\|R\|+\|B\|^2\bar c/\usig(\Sigma_W)\le\sa$, and with $\gamma^\star\ge1$ this gives $r_0\le\frac12$.

\begin{applemma}[Band of the estimated input Hessian]\label{lem:hessband}
Let Assumption~\ref{ass:pe} hold.
For every $k\ge t_0$ with $\SNR_k\ge\nu_1$ and every $K\in\cS(\bar c)$, the input Hessian $\hat S$ of $K$ on $\hat\Theta_k$ is well defined, $\|\hat S-S_K\|\le\frac{\underline{s}}8$ and $\frac{\underline{s}}2I\preceq\hat S\preceq\sa I$.
\end{applemma}

\runin{Proof of Lemma~\ref{lem:hessband}}
Write $\xi:=\SNR_k^{-1}$, which is at most $1/\nu_1$ by the hypothesis $\SNR_k\ge\nu_1$, and write $S:=S_K$, $(\hat A,\hat B):=(\hat A_k,\hat B_k)$, and $\hat P$, $\hat S$ for the objects of $K$ on $(\hat A,\hat B)$.
Since $\SNR_k\ge\nu_1\ge1/p_2$, Fact~\ref{f:pert} at index $k$ shows that $\hat A+\hat BK$ is Schur, so $\hat P$ and $\hat S$ are well defined, and that $\|\hat P-P_K\|\le p_1\xi$.
With Lemma~\ref{lem:bounds} this gives $\|\hat P\|\le\bar p$, and~\cite[Lemma~1]{zhao2025pgac} gives $\|\hat B-B\|\le\|\hat\Theta_{k}-\Theta\|\le\xi\le p_2$, so $\|\hat B\|\le\bar b$.
The difference $\hat S-S$ decomposes as
\begin{equation}\label{eq:ad1}
\begin{aligned}
\hat S-S={}&(\hat B-B)^\top\hat P\hat B+B^\top(\hat P-P_{K})\hat B\\
&+B^\top P_{K}(\hat B-B)
\end{aligned}
\end{equation}
and, since $\|\hat B-B\|\le\xi$, $\|\hat B\|\le\bar b$, $\|\hat P\|\le\bar p$, $\|\hat P-P_{K}\|\le p_1\xi$ and $\|P_{K}\|\le\bar c/\usig(\Sigma_W)$, its norm is at most $p_S\xi$, which is at most $\underline{s}/8$ because $\xi\le1/\nu_1\le\underline{s}/(8p_S)$; moreover $\|\hat S\|\le\|R\|+\bar b^2\bar p=\sa$.
With $S\succeq \underline{s}I$ this gives $\frac{\underline{s}}2I\preceq\hat S\preceq\sa I$.\hfill$\square$

\begin{applemma}[One step of Algorithm~\ref{alg:akpi}]\label{lem:onestep}
Let Assumption~\ref{ass:pe} hold, let $\alpha\in(0,1)$ satisfy (C1), let $t\ge t_0$, and suppose that $K_t\in\cS(\bar c)$, $\frac{\underline{s}}{2(1-\alpha)}I\preceq H_t\preceq\frac{\sa}{1-\alpha}I$ and $\SNR_{t+1}\ge\nu_1$.
Then $\hat S_t$ and $\hat E_t$ are well defined, and the following hold.
\pitem{(i)} $\frac{\underline{s}}2I\preceq\hat S_t\preceq\sa I$ and $\frac{\underline{s}}{2(1-\alpha)}I\preceq H_{t+1}\preceq\frac{\sa}{1-\alpha}I$.
\pitem{(ii)} $K_{t+1}\in\cS(\bar c)$ and $C(K_{t+1})-C^\star\le q_\alpha\big(C(K_t)-C^\star\big)+(1-\alpha)\nu_5\SNR_{t+1}^{-1}$.
\pitem{(iii)} $\|\Sigma_{K_{t+1}}^{-1}\Sigma_{K_t}\|\le1+\nu_4/2$.
\end{applemma}

\runin{Proof of Lemma~\ref{lem:onestep}}
The computed step is compared with the exact step $K'':=K_t-H_{t+1}^{-1}E_{K_t}$, which uses the computed preconditioner and the true slope and contracts the cost gap $C(K_t)-C^\star$ by Lemma~\ref{lem:gn}; as in~\cite{zhao2025pgac}, the difference of the two steps is bounded by a multiple of $1/\SNR_{t+1}$.
Write $\xi:=\SNR_{t+1}^{-1}$, which is at most $1/\nu_1$ by the hypothesis $\SNR_{t+1}\ge\nu_1$, and write $S:=S_{K_t}$, $E:=E_{K_t}$, $(\hat A,\hat B):=(\hat A_{t+1},\hat B_{t+1})$, and $\hat S_t$, $\hat E_t$ for the objects of $K_t$ on $(\hat A,\hat B)$.

We first prove~(i).
Since $\SNR_{t+1}\ge\nu_1\ge1/p_2$, Fact~\ref{f:pert} at index $t+1$ shows that $\hat A+\hat BK_t$ is Schur, so $\hat S_t$ and $\hat E_t$ are well defined, and that $\|\hat E_t-E\|\le p_9\xi$.
Lemma~\ref{lem:hessband} with $k=t+1$ and $K=K_t$, which applies since $K_t\in\cS(\bar c)$ and $\SNR_{t+1}\ge\nu_1$, gives $\|\hat S_t-S\|\le\frac{\underline{s}}8$ and $\frac{\underline{s}}2I\preceq\hat S_t\preceq\sa I$.
Since $\alpha\frac{\underline{s}}{2(1-\alpha)}+\frac{\underline{s}}2=\frac{\underline{s}}{2(1-\alpha)}$ and $\alpha\frac{\sa}{1-\alpha}+\sa=\frac{\sa}{1-\alpha}$, the recursion $H_{t+1}=\alpha H_t+\hat S_t$ keeps $H_{t+1}$ in the band, which proves~(i), and $M_t:=H_{t+1}^{-1}$ satisfies
\begin{equation}\label{eq:ad3}
\|M_t\|\le2(1-\alpha)/\underline{s},\qquad\lambda_{\min}(M_t)\ge(1-\alpha)/\sa.
\end{equation}

We now prove~(ii), and first show that the exact step $K''=K_t-M_tE$ contracts the cost gap $C(K_t)-C^\star$.
Since $S\preceq\hat S_t+\frac{\underline{s}}8I\preceq H_{t+1}+\frac{\underline{s}}8I$ by Lemma~\ref{lem:hessband} and $H_{t+1}=\alpha H_t+\hat S_t\succeq\hat S_t$, and $M_t^2\preceq\|M_t\|M_t$,
\begin{equation}\label{eq:ad4}
M_tSM_t\preceq M_t+\tfrac{\underline{s}}8M_t^2\preceq\big(1+\tfrac{\underline{s}}8\|M_t\|\big)M_t\preceq\tfrac32M_t.
\end{equation}
Lemmas~\ref{lem:loewner} and~\ref{lem:gn} with $K=K_t$, $M=M_t$ and $\chi=\frac32$ give $K''\in\cS$, $C(K'')\le C(K_t)$ and $C(K'')-C^\star\le(1-\frac12\varkappa/\gamma^\star)(C(K_t)-C^\star)$, and $\varkappa\ge\lambda_{\min}(M_t)\underline{s}\ge(1-\alpha)\underline{s}/\sa$ by~\eqref{eq:ad3} bounds this factor by $1-\frac{(1-\alpha)\underline{s}}{2\sa\gamma^\star}=q_\alpha$.

We then show that the computed step $K_{t+1}=K_t-M_t\hat E_t$ is a small change of the exact step $K''=K_t-M_tE$.
The two steps differ by $K_{t+1}-K''=-M_t(\hat E_t-E)$, so~\eqref{eq:ad3} and $\|\hat E_t-E\|\le p_9\xi$ give
\begin{equation}\label{eq:ad5}
\|K_{t+1}-K''\|=\|M_t(\hat E_t-E)\|\le\tfrac{2(1-\alpha)}{\underline{s}}p_9\xi\le\tfrac{p_4}{2},
\end{equation}
where the last inequality uses $\xi\le1/\nu_1\le\underline{s}p_4/(8p_9)$ and $1-\alpha\le1$.
Since $C(K'')\le C(K_t)\le\bar c$, Fact~\ref{f:lip}, applied at $K''$ with $\tilde K=K_{t+1}$, gives $K_{t+1}\in\cS$ and $|C(K_{t+1})-C(K'')|\le p_6\|K_{t+1}-K''\|\le(1-\alpha)\nu_5\xi$, by~\eqref{eq:ad5} and $\nu_5=2p_6p_9/\underline{s}$.
Combining this bound with the contraction of the exact step gives
\begin{equation}\label{eq:ad6}
C(K_{t+1})-C^\star\le q_\alpha\big(C(K_t)-C^\star\big)+(1-\alpha)\nu_5\SNR_{t+1}^{-1}.
\end{equation}
It remains to show that the cost stays below $\bar c$.
Since $\nu_5\xi\le r_0/2$ by $\xi\le1/\nu_1\le r_0/(2\nu_5)$, and $1-q_\alpha=(1-\alpha)r_0$, subtracting $C(K_t)-C^\star$ from both sides of~\eqref{eq:ad6} gives $C(K_{t+1})-C(K_t)\le(1-\alpha)r_0\big(\tfrac12-(C(K_t)-C^\star)\big)\le(1-\alpha)r_0\big(1-(C(K_t)-C^\star)\big)$.
Hence the cost does not increase if $C(K_t)-C^\star\ge1$, and otherwise $C(K_{t+1})\le C^\star+1+r_0\le\bar c$, since $r_0\le\frac12$ and $C^\star\le C(K_{t_0})$.
Thus $K_{t+1}\in\cS(\bar c)$, which together with~\eqref{eq:ad6} proves~(ii).

Finally, we prove~(iii).
By (C1), that is, $1-\alpha\le\nu_2$, and by $\|\hat E_t\|\le\|E\|+p_9\xi\le p_8$, the gain update $K_{t+1}-K_t=-M_t\hat E_t$ satisfies
\begin{equation}\label{eq:ad8}
\|K_{t+1}-K_t\|=\|M_t\hat E_t\|\le\tfrac{2(1-\alpha)p_8}{\underline{s}}\le\min\big\{p_4,\tfrac{\usig(\Sigma_W)\nu_4}{4\nu_3^2p_5}\big\},
\end{equation}
so Fact~\ref{f:lip} gives $\|\Sigma_{K_{t+1}}-\Sigma_{K_t}\|\le\usig(\Sigma_W)\nu_4/(4\nu_3^2)$, and since $\Sigma_{K_{t+1}}\succeq\usig(\Sigma_W)I_n$, $\|\Sigma_{K_{t+1}}^{-1}\Sigma_{K_t}\|\le1+\|\Sigma_{K_{t+1}}-\Sigma_{K_t}\|/\usig(\Sigma_W)$${}\le1+\nu_4/(4\nu_3^2)\le1+\nu_4/2$, since $\nu_3\ge1$.\hfill$\square$

\begin{table*}[t]
\caption{Guarantees of indirect PGAC with the Gauss--Newton step~\cite[Theorems~1 and~4]{zhao2025pgac} and of KLPI (Theorem~\ref{thm:adaptive})}
\label{tab:pgac-thm}
\centering
\renewcommand{\arraystretch}{1.15}
\begin{tabular}{@{}p{0.15\textwidth}p{0.38\textwidth}p{0.41\textwidth}@{}}
\toprule
 & Indirect PGAC, Gauss--Newton~\cite[Theorems~1 and~4]{zhao2025pgac} & KLPI, Theorem~\ref{thm:adaptive}\\
\midrule
gain update & $K_{t+1}=K_t-2\eta\,\hat S_t^{-1}\hat E_t$ & $H_{t+1}=\alpha H_t+\hat S_t$, \ $K_{t+1}=K_t-H_{t+1}^{-1}\hat E_t$\\
probing & a design parameter & $\Omega_t\sim\Nor\big(0,\tfrac{\tau_t}{2}H_t^{-1}\big)$\\
\midrule
assumed in both & \multicolumn{2}{l}{$K_{t_0}\in\cS$ and invertible $\Phi_t$ (Assumption~\ref{ass:pe})}\\
constants depend on & $(A,B,Q,R,K_{t_0})$ & $(A,B,Q,R,\Sigma_W,C(K_{t_0}))$\\
step size & $\eta\le\tilde\nu_2$ and $\eta<2\|\Sigma^\star\|$ & $\alpha<1$ and (C1) $1-\alpha\le\nu_2$\\
data condition & $\SNR_t\ge\tilde\nu_1$ for all $t$ & (C2) $\SNR_t\ge\nu_1$ for all $t\ge t_0$\\
initial preconditioner & none (no preconditioner) & (C3) $\frac{\underline{s}}{2(1-\alpha)}I\preceq H_{t_0}\preceq\frac{\sa}{1-\alpha}I$\\
\midrule
stability of the gains & sequentially stable & sequentially stable\\
state bound & $\|X_t\|\le\tilde\nu_3\big(1-\tfrac{\tilde\nu_4}{2}\big)^{t-t_0}\|X_{t_0}\|$ & $\|X_t\|\le\nu_3\big(1-\tfrac{\nu_4}{2}\big)^{t-t_0}\|X_{t_0}\|$\\
 & $\quad+\tfrac{2\tilde\nu_3}{\tilde\nu_4}\max_{t_0\le i<t}\|B\Omega_i+W_i\|$ & $\quad+\tfrac{2\nu_3}{\nu_4}\max_{t_0\le i<t}\|B\Omega_i+W_i\|$\\
contraction factor & $\rho_\eta=1-\frac{\eta}{2\|\Sigma^\star\|}$ & $q_\alpha=1-\frac{(1-\alpha)\underline{s}}{2\sa\gamma^\star}$\\
cost gap & $C(K_t)-C^\star\le\rho_\eta^{t-t_0}\big(C(K_{t_0})-C^\star\big)$ & $C(K_t)-C^\star\le q_\alpha^{t-t_0}\big(C(K_{t_0})-C^\star\big)$\\
 & $\quad+\eta\,\tilde\nu_5\sum_{i=t_0}^{t}\rho_\eta^{t-i}\,\SNR_i^{-1}$ & $\quad+(1-\alpha)\nu_5\sum_{i=t_0+1}^{t}q_\alpha^{t-i}\,\SNR_{i}^{-1}$\\
\bottomrule
\end{tabular}

\smallskip
\begin{minipage}{0.94\textwidth}\footnotesize
Constants with a tilde are those of~\cite{zhao2025pgac}, which defines its cost with $\Sigma_W=I_n$, so $C$, $C^\star$ and $\Sigma^\star$ in the PGAC column are those of $\Sigma_W=I_n$.
In both columns, $\hat S_t$ and $\hat E_t$ are the input Hessian and the slope of $K_t$ on the plant estimate $\hat\Theta_{t+1}$.
We read~\cite[Theorem~4]{zhao2025pgac} as replacing the condition $\eta\le2\|\Sigma^\star\|/\usig(R)$ of~\cite[Theorem~1]{zhao2025pgac} by $\eta<2\|\Sigma^\star\|$.
\end{minipage}
\end{table*}

\runin{Proof of Theorem~\ref{thm:adaptive}}
We show by induction on $t\ge t_0$ that $K_t\in\cS(\bar c)$ and $\frac{\underline{s}}{2(1-\alpha)}I\preceq H_t\preceq\frac{\sa}{1-\alpha}I$; this induction hypothesis holds at $t_0$ by (C3) and $C(K_{t_0})<\bar c$.
If the induction hypothesis holds at $t$, then (C2) gives $\SNR_{t+1}\ge\nu_1$, and Lemma~\ref{lem:onestep}(i) and~(ii) give the induction hypothesis at $t+1$.
Hence Lemma~\ref{lem:onestep} applies at every $t\ge t_0$, so $\hat S_t$ and $\hat E_t$ are well defined, and the bounds of Lemma~\ref{lem:onestep}(ii) and~(iii) hold at every $t\ge t_0$.
Facts~\ref{f:strong} and~\ref{f:state} then give sequential stability and~\eqref{eq:statebound}, and unrolling the cost bound of Lemma~\ref{lem:onestep}(ii) from $t_0$ gives~\eqref{eq:adrate}.
For the last claim of Theorem~\ref{thm:adaptive}, note that $K_{t_0}\in\cS(\bar c)$, since $K_{t_0}\in\cS$ and $C(K_{t_0})<\bar c$.
Lemma~\ref{lem:hessband} with $k=t_0$ and $K=K_{t_0}$, which applies by (C2), gives $\frac{\underline{s}}2I\preceq\hat S_{\mathrm{off}}\preceq\sa I$, so $H_{t_0}=\frac{\hat S_{\mathrm{off}}}{1-\alpha}$ satisfies (C3).\hfill$\square$

\runin{Comparison with PGAC}
We compare Theorem~\ref{thm:adaptive} with the guarantee of~\cite[Theorems~1 and~4]{zhao2025pgac} for PGAC's Gauss--Newton step and write both in the notation of this paper (Table~\ref{tab:pgac-thm}).
KLPI's gain update becomes this step when its preconditioner $H_{t+1}$ is replaced by $\hat S_t/(1-\alpha)$, which gives the step size $\eta=\frac{1-\alpha}2$.
At this step size both guarantees state that the gains are sequentially stable, that the state is bounded by a decaying transient plus a multiple of the largest disturbance, and that the cost gap contracts linearly up to a discounted sum of $1/\SNR_i$.
The constants of the two guarantees come from different proofs, and we do not compare them.
The two guarantees differ in four points.
First, (C3) asks the initial preconditioner $H_{t_0}$ to lie in the band in which the proof keeps $H_t$, and Theorem~\ref{thm:adaptive} shows that (C3) holds for $H_{t_0}=\hat S_{\mathrm{off}}/(1-\alpha)$ under (C2).
PGAC's Gauss--Newton step uses the current estimated input Hessian alone and has no preconditioner to initialize.
Second, the guaranteed decrease of the cost gap per sample differs by a constant factor.
When $\Sigma_W=I_n$, the covariance with which~\cite{zhao2025pgac} defines its cost, $\gamma^\star=\|\Sigma^\star\|$, and the decrease $1-q_\alpha=\frac{(1-\alpha)\underline{s}}{2\sa\|\Sigma^\star\|}$ of Theorem~\ref{thm:adaptive} is $2\underline{s}/\sa$ times the decrease $1-\rho_\eta=\frac{1-\alpha}{4\|\Sigma^\star\|}$ of~\cite[Theorem~4]{zhao2025pgac}.
Third, the sum in~\eqref{eq:adrate} starts at $i=t_0+1$, because Algorithm~\ref{alg:akpi} computes its first gain $K_{t_0+1}$ from the plant estimate $\hat\Theta_{t_0+1}$, whereas the sum in~\cite[Theorem~4]{zhao2025pgac} starts at $i=t_0$.
Finally, the condition $\eta<2\|\Sigma^\star\|$ of~\cite[Theorem~4]{zhao2025pgac} has no counterpart in Theorem~\ref{thm:adaptive}, and it holds at $\eta=\frac{1-\alpha}2$, since $\Sigma^\star\succeq I_n$ when $\Sigma_W=I_n$.

\runin{Proof of Proposition~\ref{prop:regret}}
By the Riccati equation $P^\star=Q+A^\top P^\star A-A^\top P^\star B(S^\star)^{-1}B^\top P^\star A$ and $C^\star=\Tr(P^\star\Sigma_W)$, for all $x$, $u$,
\begin{equation}\label{eq:reg1}
\begin{aligned}
&c(x,u)+\EE_{X'\sim\Nor(Ax+Bu,\Sigma_W)}\big[X'^\top P^\star X'\big]\\
&\quad=x^\top P^\star x+C^\star+(u-K^\star x)^\top S^\star(u-K^\star x).
\end{aligned}
\end{equation}
Evaluating~\eqref{eq:reg1} at $(X_t,U_t)$, which is legitimate since $W_t$ is independent of $(X_t,U_t)$, taking expectations under $\PP$, which are finite by the moment hypotheses, and summing over $t<T$ telescopes the terms in $P^\star$.
With $U_t-K^\star X_t=(K_t-K^\star)X_t+\Omega_t$, the cross term vanishes because $\Omega_t$ has zero mean given $\mathcal I_t$, which determines $(X_t,K_t,\Psi_t)$, and $\EE_{\PP}[\Omega_t^\top S^\star\Omega_t]=\EE_{\PP}[\Tr(S^\star\Psi_t)]$.\hfill$\square$

\end{adjustwidth}
\end{document}